\documentclass[a4paper,11pt]{article}
\usepackage[utf8]{inputenc}
\usepackage{enumerate}
\usepackage{lmodern}
\usepackage[T1]{fontenc}
\usepackage{verbatim}
\usepackage{textcomp}

\usepackage[english]{babel}
\usepackage[a4paper,vmargin={3.5cm,3.5cm},hmargin={2.5cm,2.5cm}]{geometry}
\usepackage[font=sf, labelfont={sf,bf}, margin=1cm]{caption}
\usepackage[pdftex]{hyperref}

\usepackage[pdftex]{color,graphicx}

\usepackage{amsmath,amsfonts,amssymb,amsthm,mathrsfs}

\newtheorem{theorem}{Theorem}
\newtheorem{definition}{Definition}
\newtheorem{lemma}[theorem]{Lemma}

\newtheorem{corollary}[theorem]{Corollary}
\newtheorem{proposition}[theorem]{Proposition}

\def\llbracket{[\hspace{-.10em} [ }
\def\rrbracket{ ] \hspace{-.10em}]}

\def\w{\mathrm{w}}

\def\t{{\mathcal T}}
\def\cc{\mathcal{C}}
\def\ve{\varepsilon}
\def\wt{\widetilde}
\def\wh{\widehat}
\def\bm{{\bf m}}

\def\bn{{\bf n}}

\def\bv{\mathbf{v}}
\def\bV{\mathbf{V}}
\def\S{\mathcal{S}}
\def\W{\mathcal{W}}

\def\mm{\mathcal{M}}
\def\dd{\mathcal{D}}

\def\la{\longrightarrow}
\def\K{{\mathbb K}}
\def\M{{\mathbb M}}

\def\E{{\mathbb E}}
\def\P{{\mathbb P}}
\def\N{{\mathbb N}}
\def\T{{\mathbb T}}

\def\R{{\mathbb R}}
\def\Q{{\mathbb Q}}
\def\F{{\mathbb F}}
\def\D{{\mathbb D}}
\def\SS{{\mathbb S}}
\def\dg{\mathrm{d}_{\mathrm{gr}}}
\def\ov{\overline}

\def\dd{\mathrm{d}}

\def\build#1_#2^#3{\mathrel{
\mathop{\kern 0pt#1}\limits_{#2}^{#3}}}

\title{Brownian geometry}
\author{Jean-Fran\c cois Le Gall\footnote{Supported by the ERC Advanced Grant 740943 {\sc GeoBrown}}}

\date{\small Universit\'e Paris-Sud}

\begin{document}
\maketitle

\begin{abstract}
We present different continuous models of random geometry that have been introduced and studied
in the recent years. In particular, we consider the Brownian map, which is the universal scaling limit
of large planar maps in the Gromov-Hausdorff sense, and the Brownian disk, which 
appears as the scaling limit of planar maps with a boundary.  We discuss the construction of these models,
and we emphasize the role played by Brownian motion indexed by
the Brownian tree.
\end{abstract}

{\small\tableofcontents}
\section{Introduction}

The goal of this work is to survey a number of recent developments 
concerning the continuous models of planar random geometry that have been
studied extensively in the last ten years, and their connections with discrete models.
A very important feature of the continuous models that we will present is their
universality, meaning that they appear in the scaling limit of many different
discrete models. This is similar of course to the universality of standard Brownian motion,
which is the scaling limit of all random walks satisfying mild moment conditions.
Partly because of this analogy, and also because Brownian motion plays
a crucial role in the construction of our basic objects of study, 
we use the name Brownian geometry for the general area 
of continuous models of random geometry that are discussed below.
In the present article, we stress the role played by Brownian 
motion indexed by the Brownian tree, which is the main ingredient of the
construction of the random metric space called the Brownian map
and of other models, and which in our opinion is also an important 
object worth of study in its own. Many properties of Brownian 
motion indexed by the Brownian tree, in particular the excursion theory 
presented in Section \ref{excu-theory} below, have direct applications
to Brownian geometry.

The discrete models of random geometry that we will consider are 
planar maps, which are finite connected graphs embedded in
the two-dimensional sphere and viewed up to direct homeomorphisms
(see Section \ref{DisGeo} below for a more precise definition). 
The faces of a planar map are the connected components of the complement of the union of edges,
and important 
particular cases of planar maps
are triangulations, respectively, quadrangulations, where all faces are bounded by $3$, resp.~$4$, edges.
We note that
many of the results that follow can be extended to
graphs embedded in surfaces of higher genus, but we will not discuss these
extensions here. Planar maps are important objects of study in combinatorices, and random planar maps have been used for a long time 
by theoretical physicists as models of random geometry, in the setting
of two-dimensional quantum gravity, see in particular \cite{Wat} and the book \cite{ADJ}.
From the mathematical point of view, a natural question is to consider a
planar map chosen at random in  a suitable class, say the class 
of all triangulations with a fixed number $n$ of faces, and to investigate the
properties of this object when $n\to \infty$. One expects, in a way similar
to the convergence of rescaled random walks to Brownian motion, 
that, when its size tends to infinity, the random planar map, suitably rescaled, 
will be close to a certain continuous model. It turns out that this vague idea can be
made precise in the framework of the Gromov-Hausdorff convergence of
compact metric spaces (see e.g.~\cite{BBI} for basic facts about the
Gromov-Hausdorff distance). Starting from a random planar map
$M_n$ uniformly distributed over the class of all triangulations with $n$
faces (or quadrangulations with $n$ faces), one shows \cite{Uniqueness,Mie2}
that the vertex set $V(M_n)$ equipped with the graph distance rescaled 
by the factor $n^{-1/4}$ converges in distribution in the 
Gromov-Hausdorff sense to a limiting random compact metric space called the 
Brownian map, see Theorem \ref{main} below (the case of
triangulations had been conjectured by Schramm \cite{Sch}). The proof of this convergence 
was strongly motivated by earlier results concerning asymptotics for the
two-point function \cite{CS} or the three-point function \cite{BG1} of random quadrangulations. The
preceding convergence to the Brownian map has
been extended to many classes of random planar maps, always with the
same limiting space, up to unimportant scaling factors on the distance: This is the
universality property of the Brownian map, which was already mentioned above.

The construction of the Brownian map, and the 
relevance of Brownian motion indexed by the Brownian tree,
are best understood from purely combinatorial considerations
about planar maps. Perhaps surprisingly, 
various classes of planar maps are in one-to-one correspondence
with certain classes of discrete trees whose vertices are
assigned integer labels. A common feature of these
bijections is the fact that labels assigned to the vertices of the tree are closely related
to graph distances from a distinguished vertex in the associated planar map.
Therefore, a good understanding of the labeled tree associated with
a random planar map yields useful information about the metric properties
of the vertex set of the planar map equipped with the graph distance.
In Section \ref{DisBij} below, we present the simplest example
of the bijections between planar maps and labeled trees, in the case of quadrangulations. 

It turns out that the tree associated with a large random planar map is
close, modulo a suitable rescaling, to the continuous random tree 
which we call the Brownian tree (this is essentially the CRT introduced and
studied by Aldous \cite{Al1,Al3}). Furthermore, labels on the tree 
behave like Brownian motion indexed by the Brownian tree  when the size of the planar map goes to infinity.
At least informally, these observations explain the construction of the Brownian map
which is presented in Section \ref{sec:consBM}: Following \cite{ALG},
we introduce the concept of a snake trajectory --- this is 
a convenient framework for studying the Brownian snake driven by a Brownian excursion \cite{Zurich},
which is basically the same object as Brownian motion indexed by the Brownian tree ---
and explain how to associate a compact metric space with a snake trajectory.
If the snake trajectory is chosen at random according to the normalized Brownian snake excursion measure, the associated compact metric space is
the Brownian map, and the ``labels'' (the values of the tree-indexed Brownian motion)
are related to distances from a distinguished point of the Brownian map.

The Brownian map is by no means the only interesting model in our Brownian geometry.
In Section \ref{Bplane}, we briefly present the Brownian plane, which is an infinite-volume version
of the Brownian map and can be obtained as the scaling limit of the infinite random lattices
called the UIPT (for uniform infinite planar triangulation) and the UIPQ 
(for uniform infinite planar quadrangulation). 
In Section \ref{PlaBou}, we introduce Brownian disks as scaling limits 
of planar maps with a boundary, when the boundary size tends to infinity \cite{Bet,BM}. 
In contrast with the Brownian map, which is homeomorphic to the two-dimensional sphere,
Brownian disks are homeomorphic to the closed unit disk. We pay special attention to the free Brownian disk,
which has a fixed boundary size or perimeter but a random volume.

Section \ref{Br-Disk} 
presents a construction of Brownian disks from a continuous random tree 
equipped with Brownian labels, which is analogous to the construction of the Brownian map,
with the difference that the labels now correspond to distances from the boundary (this is in
contrast with the previous constructions of \cite{Bet,BM}, which also used labeled trees, but with a
different interpretation of labels). Our construction relies on an excursion theory for 
Brownian motion indexed by the Brownian tree, which is developed in Section \ref{excu-theory} and is of independent 
interest. Roughly speaking,
if $\t_\zeta$ denotes the Brownian tree and $(Z_a)_{a\in\t_\zeta}$
denotes Brownian motion indexed by $\t_\zeta$, we describe the distribution of 
``excursions'' of $Z$ away from $0$, each excursion corresponding to the restriction of  $Z$ to one
connected component of $\{a\in\t_\zeta:Z_a\not =0\}$. We obtain that these excursions
are independent conditionally given their ``boundary sizes'', and distributed 
according to a certain excursion measure on snake trajectories.

The construction of Section \ref{Br-Disk} makes it possible to identify certain subsets 
of the Brownian map as Brownian disks. In particular, Theorem \ref{ccBm} shows that
connected components of the complement of the ball of radius $r$ centered at the distinguished
point in the Brownian map are independent Brownian disks, conditionally on
their boundary sizes and volumes. A similar result holds for the free Brownian disk $\D$: If $r>0$ and $H(x)$ denotes the
distance from a point $x\in\D$ to the boundary, connected components
of the set $\{x\in\D:H(x)>r\}$ are independent free Brownian disks conditionally on their boundary sizes.
Finally, in Section \ref{CutBr}, we present the very recent results of \cite{LGR} studying the sequence of boundary sizes of 
the connected components of $\{x\in\D:H(x)>r\}$  as a process parameterized by $r$. We show that
this process is a growth-fragmentation process whose distribution is completely determined. The latter
result is very closely related to the recent papers \cite{BCK,BBCK} investigating scaling limits for 
a similar process associated with triangulations with a boundary.

Even if it was not possible to provide detailed proofs in this survey, we have tried to
sketch the main ideas underlying several important results. We give a detailed
presentation of Schaeffer's bijection between quadrangulations and labeled trees, and, at the end of Section \ref{DisBij}, we
explain informally why the construction of the Brownian map, which may appear 
rather involved at first glance, is a continuous counterpart of this bijection. 
Similarly in Section \ref{Br-Disk}, we emphasize that the study of 
connected components of the complement of balls in the Brownian map
can be reduced to the study of excursions of Brownian motion
indexed by the Brownian tree. 

Let us briefly mention several recent articles that are related 
to the present work. The paper \cite{FPP} discusses the 
Gromov-Hausdorff convergence
of rescaled planar maps when the graph distance is replaced
by a ``local modification'', and shows that the scaling limit is still the Brownian map.
The study of the UIPT and the UIPQ has given rise to a number
of interesting developments: See in particular \cite{GN} for a proof of the recurrence 
of simple random walk on these infinite random lattices.
Hyperbolic versions of the Brownian plane have been studied
by Budzinski \cite{Bud}. The Brownian half-plane, which also appears as the scaling limit of quadrangulations 
with a boundary when the volume and the boundary size tend to infinity in a suitable way, is discussed in \cite{GM1} and \cite{BMR} 
--- a presumably equivalent construction had been given earlier by Caraceni and Curien \cite{CC}. 
 The paper \cite{BMR} provides an exhaustive study of possible 
scaling limits of quadrangulations with a boundary, leading to new models of Brownian geometry in addition to
the Brownian map, the Brownian plane or the Brownian disk. In a series of recent papers,
Miller and Sheffield \cite{MS0,MS1,MS2,MS3} have developed a completely new approach 
to the Brownian map, showing also that this random compact metric space can be equipped with a conformal structure which 
is linked to Liouville quantum gravity. An important step in this approach \cite{MS0} was the derivation of
an axiomatic characterization of the Brownian map. The paper \cite{MS0} uses a definition of Brownian disks 
which is different from the one in \cite{Bet} but which can be shown to be consistent with the latter thanks
to the results of \cite{Disk}. Brownian disks play an important role in the recent work
\cite{GM0,GM1,GM3} of Gwynne and Miller motivated by the study of 
statistical physics models on random planar maps. We also mention 
the paper \cite{GMS} showing that certain discrete {\it conformal} embeddings 
of random planar maps converge to their continuous counterparts. Finally we refer 
to \cite{LGM2,Mie-Saint-Flour} for pedagogical presentations of random planar maps and
the convergence to the Brownian map. 

\smallskip
\noindent {\it Acknowledgement.} It is a pleasure to thank the organizers of the Takagi lectures for 
giving me the opportunity to discuss the present work at this prestigious meeting.

\section{Discrete and continuous models of random geometry}
\label{DisGeo}

\subsection{Planar maps}

The basic discrete model of random geometry that we will consider is a 
random planar map. Let us start with a precise definition.

\begin{definition}
\label{planar-map}
A planar map is a proper embedding of a finite connected graph in the
two-dimensional sphere $\SS^2$. Two planar maps are identified if they
correspond via an orientation-preserving homeomorphism of the sphere.
\end{definition}

In this definition,  ``proper'' means that edges are not allowed to cross. 
The identification modulo homeomorphisms is interpreted by saying that we are only interested in the shape of the
embedding, and not in its precise details.

\begin{figure}[!h]
 \begin{center}
 \includegraphics[width=7cm]{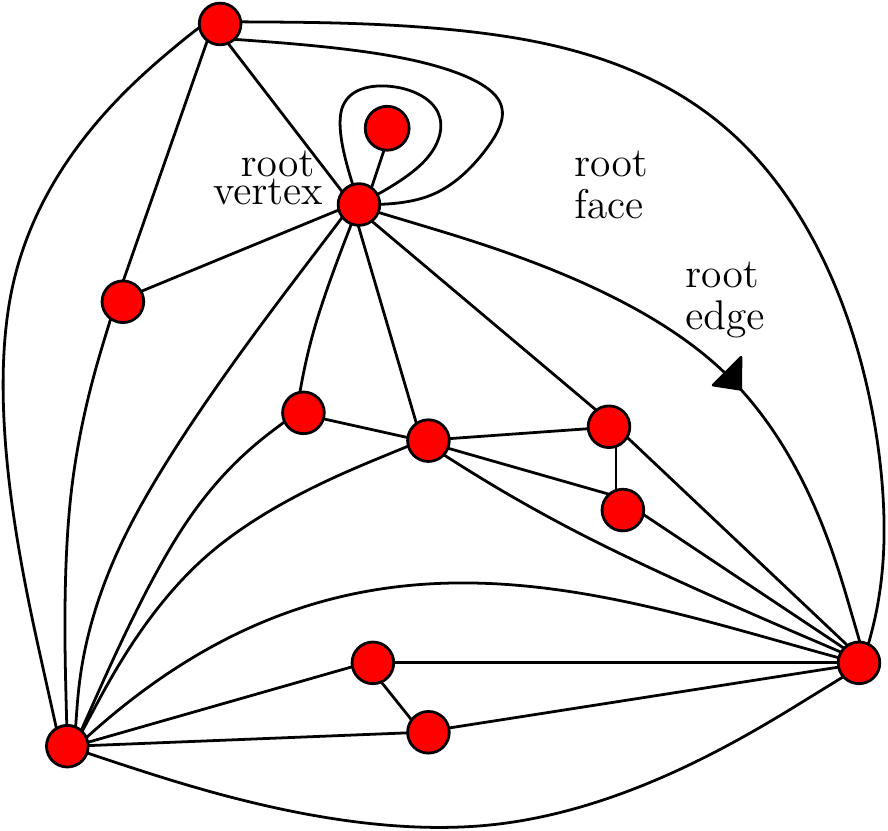}
 \caption{\label{triangu}
 A rooted triangulation with $20$ faces.}
 \end{center}
 \vspace{-5mm}
 \end{figure}

In the preceding definition, we should in fact have written ``multigraph'' instead
of graph, meaning that we allow self-loops and multiple edges. Many of the results
that follow are expected to hold, and sometimes have been proved, also for
simple planar maps where self-loops and multiple edges are forbidden, but the technicalities 
become more difficult.
See Fig.~\ref{triangu} for an example with a self-loop and a double edge.

 Thanks to the fact that the graph is embedded, we can define the notion of a face.
 Faces are the connected components of the complement of edges, or equivalently
 the regions bounded by the edges. The degree of a face is the number
 of half-edges incident to this face: Note that we say half-edges instead of edges because if
both sides of an edge are incident to the same face, this edge is counted twice in the degree (for instance the face 
inside the self-loop in Fig.~\ref{triangu} has degree $3$ though there are only two edges in its boundary). 

If $p\geq 3$ is an integer, a planar map is called a {\it $p$-angulation} if all its faces have degree $p$,
and we say {\it triangulation} when $p=3$, {\it quadrangulation} when $p=4$. Fig.~\ref{triangu} shows a triangulation
with $20$ faces.

We will deal with {\it rooted} planar maps, meaning that we distinguish an oriented edge, which is 
called the root edge. The origin of the root edge is called the root vertex, and the face lying to the left of
the root edge (this makes sense because the root edge is oriented) is called the root face. See again Fig.~\ref{triangu}.
Notice that in order to identify two rooted planar maps via an orientation-preserving homeomorphism
we require that this homeomorphism preserves the root edge. The reason for dealing with rooted maps
comes from the fact that enumeration questions, or bijections between maps and simpler objects such as trees, become more tractable
(rooting a map avoids problems related to the presence of symmetries). However, it is strongly believed that the 
results that follow hold as well for planar maps that are not rooted.

Let $p\geq 3$ and $n\geq 1$ be integers. The set of all rooted $p$-angulations with $n$ faces will be denoted
by $\M^p_n$. It is easy to see that $\M^p_n$ is empty if $p$ and $n$ are both odd integers. So when
$p$ is odd, in particular when $p=3$, we will implicitly restrict our attention to even values of $n$.
Thanks to the identification in Definition \ref{planar-map}, the set $\M^p_n$ is finite, and so it makes sense
to choose a rooted $p$-angulation with $n$ faces uniformly at random.

If $M$ is a planar map, we will denote the vertex set of $M$ by $V(M)$. We equip $V(M)$
with the usual graph distance $\dg^M$: If $v$ and $v'$ are two vertices of $M$, $\dg^M(v,v')$
is the minimal number of edges on a path from $v$ to $v'$. Our first goal is to study the metric
space $(V(M),\dg^M)$ when $M$ is chosen uniformly at random in $\M^p_n$ (for some fixed $p$)
and when $n$ is large. For this study, we will need a notion of convergence of 
a sequence of compact metric spaces.

\subsection{The Gromov-Hausdorff distance}

Let us first recall that, 
if $K_1$, $K_2$ are two compact subsets of a metric space $(E,\dd)$, the
Hausdorff distance between $K_1$ and $K_2$ is defined by
$$\dd^E_{\rm Haus}(K_1,K_2)=\inf\{\ve>0: K_1\subset U_\ve(K_2)
\hbox{ and } K_2\subset U_\ve(K_1)\}$$
where $U_\ve(K_1)=\{x\in E:\dd(x,K_1)\leq \ve\}$ is the $\ve$-enlargement of $K_1$.

\begin{definition}[Gromov-Hausdorff distance]
\label{Gro}
Let $(E_1,\dd_1)$ and $(E_2,\dd_2)$ be two compact metric spaces. The Gromov-Hausdorff distance
between $E_1$ and $E_2$ is
$$\dd_{\rm GH}(E_1,E_2)=\inf\{\dd^E_{\rm Haus}(\psi_1(E_1),\psi_2(E_2))\}$$
where the infimum is over all  isometric embeddings $\psi_1:E_1\to E$ and $\psi_2:E_2\to E$
of $E_1$ and $E_2$ into the same metric space $(E,\dd)$.
\end{definition}

Let $\K$ stand for the set of all compact metric spaces, where as usual two compact metric
spaces are identified if they are isometric. Then the Gromov-Hausdorff distance $\dd_{\rm GH}$
is a metric on $\K$, and furthermore $(\K,\dd_{\rm GH})$ is complete and separable.
In other words, $(\K,\dd_{\rm GH})$ is a Polish space, which makes it especially
suitable to study the convergence in distribution of random variables with values in $\K$. 

One can prove \cite{Mie-Perso} that a sequence $(E_n)$ of compact metric spaces converges to
a limiting space $E_\infty$ in $\K$ if and only if all spaces $E_n$ and the limit $E_\infty$
can be embedded isometrically in the same metric space $E$ in such a way that
the convergence holds in the sense of the Hausdorff distance.

\subsection{Convergence to the Brownian map}

We will now discuss the convergence in distribution of 
$(M_n,n^{-1/4}\dg^{M_n})$ when $M_n$ is chosen uniformly at random in $\M^p_n$ (for some fixed $p$).
Note that we rescale the graph distance $\dg^{M_n}$ by the factor $n^{-1/4}$: The need for such a rescaling
is clear since one expects that the diameter of the graph blows up when the number of faces
grows to infinity. The reason why the correct rescaling factor is $n^{-1/4}$ is 
more mysterious and will be best understood from the bijections between planar maps and
labeled trees that are described below (see the beginning of Section \ref{proof-main}). 

The following theorem is proved in \cite{Uniqueness}. The particular case 
of quadrangulations $p=4$ was obtained independently by Miermont \cite{Mie2}.
The case $p=3$ solves a problem of Schramm \cite{Sch}.

\begin{theorem}[The scaling limit of $p$-angulations]
\label{main}
Suppose that either $p=3$ (triangulations) or $p\geq 4$ is even. Set 
$$c_3= 6^{1/4}\quad,\quad c_p= \Big(\frac{9}{p(p-2)}\Big)^{1/4}\quad\hbox{if }p\hbox{ is even.}$$
For every integer $n\geq 2$ ($n$ even if $p=3$), let $M_n$
be uniformly distributed over $\M^p_n$. Then,
$$(V(M_n), c_p\,{n^{-1/4}}\,\dg^{M_n})
\build{\longrightarrow}_{n\to\infty}^{\rm (d)} ({\bf m}_\infty,\dd_\infty)$$
in the Gromov-Hausdorff sense. The limit
$({\bf m}_\infty,\dd_\infty)$ is a random compact metric space (that is, a random variable with values 
in $\K$) that
does not depend on $p$ and is 
called the Brownian map.
\end{theorem}

The name ``Brownian map'' is due to Marckert and Mokkadem \cite{MM1} who obtained a weak
form of the theorem in the case of quadrangulations.
We note that the role of the constants $c_p$ is only to ensure that the limit does not depend on $p$.
It is expected that the result of the theorem holds for all values of $p\geq 3$, but the
case of odd values $p\geq 5$ seems more difficult to handle for technical reasons.

The fact that the limit does not depend on $p$ is a very important feature of Theorem \ref{main}.
Roughly speaking, it means that in large scales the metric properties of a typical (large) planar map
are the same if this planar map is a triangulation, or a quadrangulation, or a $p$-angulation. 
This is the universality property of the Brownian map, which has been confirmed
in many subsequent works: In particular, analogs of Theorem \ref{main}, always with the
same limit $({\bf m}_\infty,\dd_\infty)$ hold for general planar maps with a fixed number of
edges \cite{BJM}, for bipartite planar maps with a fixed number of
edges \cite{Abr}, for simple triangulations or quadrangulations (where self-loops and multiple edges
are not allowed) \cite{AA}, for planar maps with a prescribed degree sequence \cite{Mar}, etc.
We also mention that results similar to Theorem \ref{main} hold if the graph distance
is replaced by a ``local modification'': The paper \cite{FPP} considers the so-called 
first-passage percolation distance on random triangulations 
(independent random weights are assigned to
the edges and the distance between two vertices is the minimal total weight 
of a path between them). Perhaps surprisingly, this local modification does
not change the scaling limit, which is still the Brownian map up to
a deterministic scale factor for the distance.

As a general principle, the scaling limit of large random planar maps is expected to be the Brownian map whenever some
bound is assumed on the degree of faces. On the contrary, if one considers probability distributions
on planar maps that favor the appearance of very large faces, different scaling limits may occur
(the so-called stable maps of \cite{LGM}), but we will not discuss this case here.

It is implicit in Theorem \ref{main} that the limit $({\bf m}_\infty,\dd_\infty)$ is not the degenerate space
with a single point. We make this more explicit in the following two theorems that give some
useful information about the Brownian map.

\begin{theorem} [\cite{Invent}]
\label{Hausd}
The Hausdorff dimension of $({\bf m}_\infty,\dd_\infty)$ is a.s. equal to $4$.
\end{theorem}

\begin{theorem} [\cite{LGP}]
\label{topo}
The compact metric space $({\bf m}_\infty,\dd_\infty)$ is a.s. homeomorphic to the $2$-sphere $\SS^2$.
\end{theorem}

Both these theorems can be deduced from the construction of the Brownian map from Brownian motion indexed 
by the Brownian tree that will
be given below. The proof of Theorem \ref{Hausd} is in fact relatively easy, but that of Theorem \ref{topo}
is more intricate and relies in part on an old theorem of Moore giving conditions for a quotient space
of the sphere to be homeomorphic to the sphere.

Since planar maps are defined as graphs embedded in the sphere, and since we take a limit 
where the number of vertices tends to infinity, it is maybe not surprising that the limiting metric space
has the topology of the sphere. Still, Theorem \ref{topo} implies a non-trivial combinatorial fact
about the non-existence of small  ``bottlenecks'' in a large planar map: Informally, for a random triangulation with $n$
faces, the probability that there exists a cycle with length $o(n^{1/4})$ such that both sides
of the cycle (meaning both components of its complement) have a diameter greater than $\delta n^{1/4}$,
for some fixed $\delta>0$, will tend to $0$ as $n\to\infty$. The question of the existence of small
separating cycles in random planar maps has been investigated recently in connection 
with isoperimetric inequalities \cite{LGL}.

\section{The construction of the Brownian map}
\label{sec:consBM}

In this section, we present a construction of the limiting space $({\bf m}_\infty,\dd_\infty)$ of
Theorem \ref{main}. This construction relies on the notion of Brownian motion 
indexed by the Brownian tree. We start by a brief presentation of the Brownian tree.

\subsection{The Brownian tree}
\label{sub:BrT}

Recall that an $\R$-tree is a metric space $(\t,\dd)$ such that, for every $a,b\in\t$ there is, up
to reparameterization, a 
unique continuous injective path $\gamma$ from $a$ to $b$, and the range of $\gamma$,
which will be denoted by $\llbracket a,b\rrbracket$,
is isometric to the  line segment $[0,\dd(a,b)]$. An $\R$-tree $\t$ is rooted if there is
a distinguished point $\rho\in\t$, which is called the root. This makes it possible to
define a notion of genealogy in the tree $\t$: If $a,b\in\t$, we say that $b$ is a
descendant of $a$, or $a$ is an ancestor of $b$, if $a\in \llbracket \rho,b\rrbracket$. 

In the present work, we will consider only
compact $\R$-trees, and we will use the fact that such trees can be coded by
continuous functions. Let $h:\R_+\to\R_+$ be a nonnegative continuous function on $\R_+$ such that $h(0)=0$. We assume that
$h$ has compact support, so that
$$\sigma_h:=\sup\{t\geq 0: h(t)>0\} <\infty.$$
Here and later we make the convention that $\sup\varnothing = 0$. 

For every $s,t\in\R_+$, we set
$$\dd_h(s,t):=h(s)+h(t)-2\min_{s\wedge t\leq r\leq s\vee t} h(r).$$
We note that $\dd_h$ is a pseudo-metric on $\R_+$, and thus we may introduce the associated
 equivalence relation on $\R_+$, defined by setting
$s\sim_h t$ if and only if $\dd_h(s,t)=0$, or equivalently
$$h(s)=h(t)=\min_{s\wedge t\leq r\leq s\vee t} h(r).$$
Then, $\dd_h$ induces a distance on the quotient space $\R_+/\!\sim_h$.

\begin{lemma}
\label{cod-tree}
{\rm\cite{DLG}}
The quotient space $\t_h:=\R_+/\!\sim_h$ equipped with the distance $\dd_h$ is a compact $\R$-tree called the tree 
coded by $h$. The canonical projection from $\R_+$ onto $\t_h$
is denoted by $p_h$. By definition, $\t_h$ is rooted at $\rho=p_h(0)$. 
\end{lemma}

\noindent{\bf Remark.} It is not hard to verify that any compact $\R$-tree can be represented as $\t_h$
for some (not unique) function $h$, but we will not need this fact. 

\medskip

It is often convenient to equip $\t_h$ with a {\it volume measure}, which is defined as the
push forward of Lebesgue measure on $[0,\sigma_h]$ under $p_h$.

The coding by a function makes it possible to define ``lexicographical'' intervals
on the tree. Let us explain this. If $s,t\geq 0$ and $s>t$, we make the convention that $[s,t]=[s,\infty)\cup [0,t]$
(of course, if $s\leq t$, $[s,t]$ is the usual interval). If $a,b\in \t_h$, there is a smallest
``interval'' $[s,t]$ with $s,t\geq 0$ (but not necessarily $s\leq t$) such that 
$p_h(s)=a$ and $p_h(b)=t$, and we then 
set $[a,b]=p_h([s,t])$. Note that $[a,b]$ is typically different from $[b,a]$. Intuitively, $[a,b]$ is the set
of all points of $\t_h$ that are visited when going from $a$ to $b$ around the tree
in ``clockwise order''.

Let us now randomize $h$. We let $\mathbf{n}(\mathrm{d}h)$ stand for It\^o's excursion measure
of positive excursions of linear Brownian motion (see e.g. \cite[Chapter XII]{RY}) normalized so
that, for every $\ve >0$,
$$\mathbf{n}\Big(\max_{s\geq 0} h(s) >\ve\Big) = \frac{1}{2\ve}.$$ 
Under $\bn(\dd h)$, we will write $\sigma=\sigma_h$ for the duration of the excursion $h$.
It will also be convenient to introduce the conditional probability measure $\bn_{(s)}:=\bn(\cdot\mid \sigma=s)$,
for every $s>0$. In particular $\bn_{(1)}$ is the law of the normalized excursion, and we have
$$\bn = \int_0^\infty \bn_{(s)}\, \frac{\dd s}{2\sqrt{2\pi s^3}}.$$
 
 \begin{definition}
 \label{BrT}
 The Brownian tree is the tree $\t_h$ coded by $h$ under $\bn(\dd h)$.
 \end{definition}
 
 It is important to realize that $\bn$ is an infinite measure. We can also consider 
 the tree $\t_h$ under the probability measure $\bn_{(1)}(\dd h)$, and this random tree is Aldous' continuum random tree,
 also called the CRT
 (our normalization is slightly different from the one in \cite{Al1,Al3}). However, it is often more convenient to argue under the infinite measure $\bn$. 
 
 \subsection{Snake trajectories}
 \label{sec:sna}
 
 We now propose to discuss Brownian motion indexed by the Brownian tree of 
 Definition \ref{BrT}. The fact that we are interested in a random process indexed by
 a random set creates some technical difficulties, which we will avoid here 
 by introducing the concept of a {\it snake trajectory}. 
 
A finite real path is a continuous mapping $\w:[0,\zeta_{(\w)}]\la \R$, where the 
number $\zeta_{(\w)}\geq 0$ is called the lifetime of $\w$.  We let 
$\W$ denote the space of all finite paths in $\R$. The set $\W$ is a Polish space when equipped with the
distance
$$\dd_\W(\w,\w')=|\zeta_{(\w)}-\zeta_{(\w')}|+\sup_{t\geq 0}|\w(t\wedge
\zeta_{(\w)})-\w'(t\wedge\zeta_{(\w')})|.$$
The endpoint or tip of the path $\w$ is denoted by $\wh \w=\w(\zeta_{(\w)})$.
For every $x\in\R$, we set $\W_x=\{\w\in\W:\w(0)=x\}$. The trivial element of $\W_x$ 
with zero lifetime is identified with the point $x$ --- in this way we view $\R$
as the subset of $\W$ consisting of all finite paths with zero lifetime. 
%We will also use the notation
%$\underline\w=\min\{\w(t):0\leq t\leq \zeta_{(\w)}\}$.

%We next turn to snake trajectories.

\begin{definition}
\label{def:snakepaths}
Let $x\in \R$. A snake trajectory with initial point $x$ is a continuous mapping
\begin{align*}
\omega:\ &\R_+\to \W_x\\
&s\mapsto \omega_s
\end{align*}
which satisfies the following two properties:
\begin{enumerate}
\item[\rm(i)] We have $\omega_0=x$ and the number $\sigma(\omega):=\sup\{s\geq 0: \omega_s\not =x\}$,
called the duration of the snake trajectory $\omega$,
is finite. 
\item[\rm(ii)] For every $0\leq s\leq s'$, we have
$$\omega_s(t)=\omega_{s'}(t)\;,\quad\hbox{for every } 0\leq t\leq \min_{s\leq r\leq s'} \zeta_{(\omega_r)}.$$
\end{enumerate} 
\end{definition}

Property (i) implies in particular that the function $s\mapsto\zeta_{(\omega_s)}$ has compact support.

\smallskip

\noindent{\bf Important remark.} A snake trajectory $\omega$ is completely determined by
the knowledge of the lifetime function $s\mapsto \zeta_{(\omega_s)}$
and the tip function $s\mapsto \wh \omega_s=\omega_s(\zeta_{(\omega_s)})$. Indeed,
for any $s\geq 0$ and $r\in[0,\zeta_{(\omega_s)}]$, if $\theta_s(r)=\inf\{u\geq s:\zeta_{(\omega_u)}=r\}$,
property (ii) implies that $\omega_s(r)=\wh \omega_{\theta_s(r)}$. 

\smallskip

We write $\S_x$ for the set of all snake trajectories  with initial point $x$, and 
$$\S:=\bigcup_{x\in \R} \S_x$$
for the set of all snake trajectories. 

 Let $\omega\in\S$. Then the real function 
 $s\mapsto \zeta_{(\omega_s)}$ satisfies the conditions required to define the 
 tree coded by this function (cf. Section \ref{sub:BrT}) and we will write $\t_\zeta$
 for this tree, and $p_\zeta$ for the canonical projection from $\R_+$ onto $\t_\zeta$. We sometimes say that
 $\t_\zeta$ is the genealogical tree of the snake trajectory $\omega$. Property (ii) in Definition \ref{def:snakepaths} 
 implies that $\omega_s=\omega_{s'}$ whenever $p_\zeta(s)=p_\zeta(s')$. In other words, $\omega_s$
only depends on the equivalence class of $s$ in the quotient space $\t_{\zeta}$, and the 
mapping $s\mapsto \omega_s$ induces a function defined on the 
genealogical tree $\t_{\zeta}$. We should think of the collection
$(\omega_s)_{s\geq 0}$ as forming a ``tree of paths'' whose genealogy is prescribed by $\t_\zeta$ (see
the left side of Fig.~\ref{Exc} below for an illustration).

\smallskip
\noindent{\bf Notation.} In what follows, we will consider snake trajectories
$\omega$ that may be deterministic or chosen according to a measure on $\S$, and we will use the notation 
$W_s=W_s(\omega)=\omega_s$, and $\zeta_s=\zeta_s(\omega)=\zeta_{(\omega_s)}$.

\subsection{Constructing a compact metric space from a snake trajectory}
\label{subsec:consBM}

The Brownian map of Theorem \ref{main} is constructed 
from a random snake trajectory distributed according to 
a certain probability measure.
To explain this construction, it is best to consider first the case
of a deterministic snake trajectory $\omega$. 

So we fix $\omega\in\S_0$ and we recall that
$\t_\zeta$ is the tree coded by $(\zeta_s)_{s\geq 0}$
(we use the notation explained at the end of Section
\ref{sec:sna}). If $a\in \t_\zeta$,
we set $Z_a=\wh W_s$ if $s$ is such that $p_\zeta(s)=a$
and we also say that
$W_s$ is the historical path of $a$ (by preceding
observations, this does not depend on the choice of $s$
such that $p_\zeta(s)=a$). We view 
$(Z_a)_{a\in\t_\zeta}$ as a collection of labels assigned 
to the points of $\t_\zeta$. Note that the function $a\mapsto Z_a$ is continuous on $\t_\zeta$.

We will now associate a metric space with the space trajectory $\omega$,
and roughly speaking this metric space will be obtained from
the genealogical tree $\t_\zeta$ by gluing together certain pairs of points.
Let us turn to a precise definition. For every $a,b\in\t_\zeta$, we set
\begin{equation}
\label{Dzero}
D^\circ(a,b)=Z_a + Z_b -2\max\Big(\min_{c\in[a,b]} Z_c, \min_{c\in[b,a]} Z_c\Big),
\end{equation}
where we recall that $[a,b]$ stands for the lexicographical interval from 
$a$ to $b$ in $\t_\zeta$. We note that $D^\circ(a,b)=0$ if and only if
\begin{equation}
\label{identi-BM}
Z_a = Z_b = \max\Big(\min_{c\in[a,b]} Z_c, \min_{c\in[b,a]} Z_c\Big),
\end{equation}
which informally means that $a$ and $b$ have the same label and that
we can go from $a$ to $b$ around the tree (clockwise or counterclockwise)
visiting only points whose label is at least as large as the label of $a$ and $b$. 
We then let $D(a,b)$ be the largest symmetric function of the pair $(a,b)$ that
is bounded above by $D^\circ(a,b)$ and satisfies the triangle inequality: For every
$a,b\in\t_\zeta$, 
\begin{equation}
\label{formulaD}
D(a,b) = \inf\Big\{ \sum_{i=1}^k D^\circ(a_{i-1},a_i)\Big\},
\end{equation}
where the infimum is over all choices of the integer $k\geq 1$ and of the
elements $a_0,a_1,\ldots,a_k$ of $\t_\zeta$ such that $a_0=a$
and $a_k=b$. We note that 
\begin{equation}
\label{lowerB-D}
D(a,b)\geq |Z_a-Z_b|
\end{equation}
as an immediate consequence of the similar bound for $D^\circ$.

We now observe that $D$ is a pseudo-metric on $\t_\zeta$, and 
we let $\mm$ be the associated quotient space, which is  the quotient of
$\t_\zeta$ for the equivalence relation $a\approx b$ if and only if $D(a,b)=0$.
We equip $\mm$
with the distance induced by $D$, for which we keep the same notation.
We note that $(\mm,D)$ is a compact metric space, and we
let $\Pi$ denote the canonical projection from $\t_\zeta$ onto $\mm$.
By abuse of notation, for every $x\in \mm$, we write $Z_x=Z_a$
if $x=\Pi(a)$ (by \eqref{lowerB-D} this does not depend on the choice of $a$
such that $x=\Pi(a)$). So labels can also be viewed as defined on the
quotient space $\mm$. Later it will be convenient to have a volume 
measure $\mathbf{v}(\dd x)$ on $\mm$, which is defined as the push forward of the volume 
measure on $\t_\zeta$ under $\Pi$.

The preceding construction obviously depends on the choice of $\omega$, which was
fixed in the beginning of this section. We
claim that it does so in a measurable way.

\begin{lemma}
\label{measura}
The mapping $\omega\mapsto (\mm, D)$ defined above, with values in the
space $(\K,\dd_{GH})$, is measurable.
\end{lemma}

We refer to \cite[Lemma 6]{Disk} for the proof of a more precise statement. 

Let us mention some properties of $D$ that will play a role later. Let $a_*$
be any point of $\t_\zeta$ such that
$$Z_{a_*}=\inf_{a\in \t_\zeta} Z_a.$$
The existence of such a point follows from a compactness argument (notice that $a_*$ may not be unique, but
if will follow from \eqref{dist-min} below that $\Pi(a_*)$ is uniquely determined). Then 
we have, for every $a\in\t_\zeta$,
$$D(a_*,a) = Z_a - Z_{a_*}.
$$
The lower bound $D(a_*,a) \geq Z_a - Z_{a_*}$ is immediate from \eqref{lowerB-D}.
The corresponding upper bound is also trivial since it is
clear that $D^\circ(a_*,a)= Z_a - Z_{a_*}$. So setting $Z_*=Z_{a_*}$ and $x_*=\Pi(a_*)$,
we get that, for every $x\in\mm$,
\begin{equation}
\label{dist-min}
D(x_*,x)= Z_x -Z_*.
\end{equation}
We interpret this by saying that 
$\mm$ has a distinguished point $x_*$
such that labels exactly correspond to distances from $x_*$,
up to the shift by $Z_*$. 

\subsection{Measures on snake trajectories}

We start with a key lemma.

\begin{lemma}
\label{snake-driven}
Let $h:\R_+\to \R_+$ be a continuous function with compact support such that $h(0)=0$. 
Assume that $h$ is H\"older continuous, meaning that there exist positive constants $\delta\in(0,1]$ and $C$
such that $|h(s)-h(s')|\leq C|s-s'|^\delta$ for every $s,s'\geq 0$. Then there exists a random snake trajectory
$W^h=(W^h_s)_{s\geq 0}$ with initial point $0$ such that:
\begin{itemize}
\item[\rm(i)] $\zeta_{(W^h_s)}=h(s)$, for every $s\geq 0$, a.s.
\item[\rm (ii)] The process $(\wh W^h_s)_{s\geq 0}$ is a centered Gaussian process
with covariance
$${\rm cov}(\wh W^h_s,\wh W^h_{s'})=\min_{s\wedge s'\leq r\leq s\vee s'} h(r).$$
\end{itemize}
The process $(W^h_s)_{s\geq 0}$ is called the Brownian snake driven by the function $h$. 
\end{lemma}

We note that
the distribution of $W^h$ is completely determined by properties (i) and (ii), thanks to the remark
following Definition \ref{def:snakepaths}. The intuition underlying the preceding definition
is as follows: For every $s\geq 0$, $W^h_s$ is a Brownian path with lifetime $h(s)$, when 
$h(s)$ decreases the path $W^h_s$ is erased from its tip and when $h(s)$ increases
the path $W^h_s$ is extended by adding ``little pieces of Brownian paths'' at its tip.

The proof of Lemma \ref{snake-driven} is straightforward. Note that the H\"older
continuity assumption of $h$ is used to warrant the existence of a 
continuous modification of a process satisfying properties (i) and (ii)
of Lemma \ref{snake-driven}.

As a consequence of (ii), we have $W^h_0=0$ and
$$\E[ (\wh W^h_s-\wh W^h_{s'})^2] = \dd_h(s,s'),$$
where the pseudo-metric $\dd_h$ was defined in Section \ref{sub:BrT}. 
Since we already noted that the snake trajectory $W^h$ can be viewed as
indexed by the tree $\t_h$, the last display justifies the fact that $\wh W^h$
is interpreted as Brownian motion indexed by $\t_h$. In fact, if $\varphi:[0,u]\la \t_h$
is an isometry mapping the interval $[0,u]$ onto a line segment of $\t_h$, we immediately see that
$(\wh W^h_{\varphi(r)}-\wh W^h_{\varphi(0)})_{0\leq r\leq u}$ is a linear Brownian motion. 

 If $h$ satisfies the properties in Lemma \ref{snake-driven},
we let $\mathbf{P}^h$ stand for the distribution of $W^h$, which is thus
a probability measure on the space $\S_0$.
We now introduce Brownian snake excursion measures, which will play a major
role in what follows.

\begin{definition}
\label{BSexc}
The Brownian snake excursion measure $\N_0$
is the $\sigma$-finite measure on $\S_0$
defined by
$$\N_0(\dd \omega)=\int \bn(\dd h)\, \mathbf{P}^h(\dd\omega).$$
Similarly, the normalized Brownian snake excursion measure
is the probability measure on $\S_0$ defined by
$$\N^{(1)}_0(\dd \omega)=\int \bn_{(1)}(\dd h)\, \mathbf{P}^h(\dd\omega).$$
\end{definition}

In other words, to construct a random snake trajectory distributed according 
to $\N_0$ (resp. according to $\N_0^{(1)}$) we just pick a 
Brownian excursion $h$ distributed according to $\bn$ (resp. a normalized Brownian
excursion) and consider the Brownian snake driven by $h$. This makes sense because we
know that $h$ is H\"older continuous, $\bn(\dd h)$ or $\bn_{(1)}(\dd h)$ a.e.

Lemma \ref{measura} now allows us to set the following definition.

\begin{definition}
\label{def:BroMap}
The Brownian map is the random compact metric space $(\mm,D)$
obtained via the construction of Section \ref{subsec:consBM} from a
snake trajectory $\omega$ distributed according to $\N^{(1)}_0$.
\end{definition}

One can prove \cite{Invent} that $\N^{(1)}_0$ a.s., for every $a,b\in\t_\zeta$ the property 
$D(a,b)=0$ holds if and only if $D^\circ(a,b)=0$ (the fact that
$D^\circ(a,b)=0$ implies $D(a,b)=0$ is obvious since $D\leq D^\circ$). So the 
construction of the Brownian map can be summarized by saying that
we start from the CRT $\t_\zeta$ equipped with ``Brownian labels''
$(Z_a)_{a\in\t_\zeta}$, and we identify points $a$ and $b$ of the CRT if
and only if $D^\circ(a,b)=0$, which has a simple interpretation as 
explained above after \eqref{identi-BM} (furthermore the metric $D$ is the largest metric 
bounded above by $D^\circ$). 

It is often useful to consider also the {\it free Brownian map}, which
is just the metric space $(\mm,D)$ under the measure $\N_0$. Many properties
of the free Brownian map are ``nicer'' than those of the ``standard 
Brownian map'' because there is no constraint
on the total volume, but the price to pay is to work under a 
$\sigma$-finite measure.

\section{Discrete bijections with trees}
\label{DisBij}

\subsection{Schaeffer's bijection}
\label{Scha-bi}

In this section, we explain a bijection between quadrangulations and (discrete) labeled trees, which
can be found in \cite{CS} and is
in some sense a discrete analog of the construction of the Brownian map that was
given in the previous section. In fact this discrete bijection (and its generalizations) 
plays a major role in the proof of Theorem \ref{main}, and helps to understand the
definition of the Brownian map and of the metric $D$. We restrict our attention to the case
of quadrangulations because the description is simpler in that case, but we immediately mention that similar bijections exist for
more general planar maps (see in particular \cite{BDG}).

\begin{figure}[!h]
 \begin{center}
 \includegraphics[width=15cm]{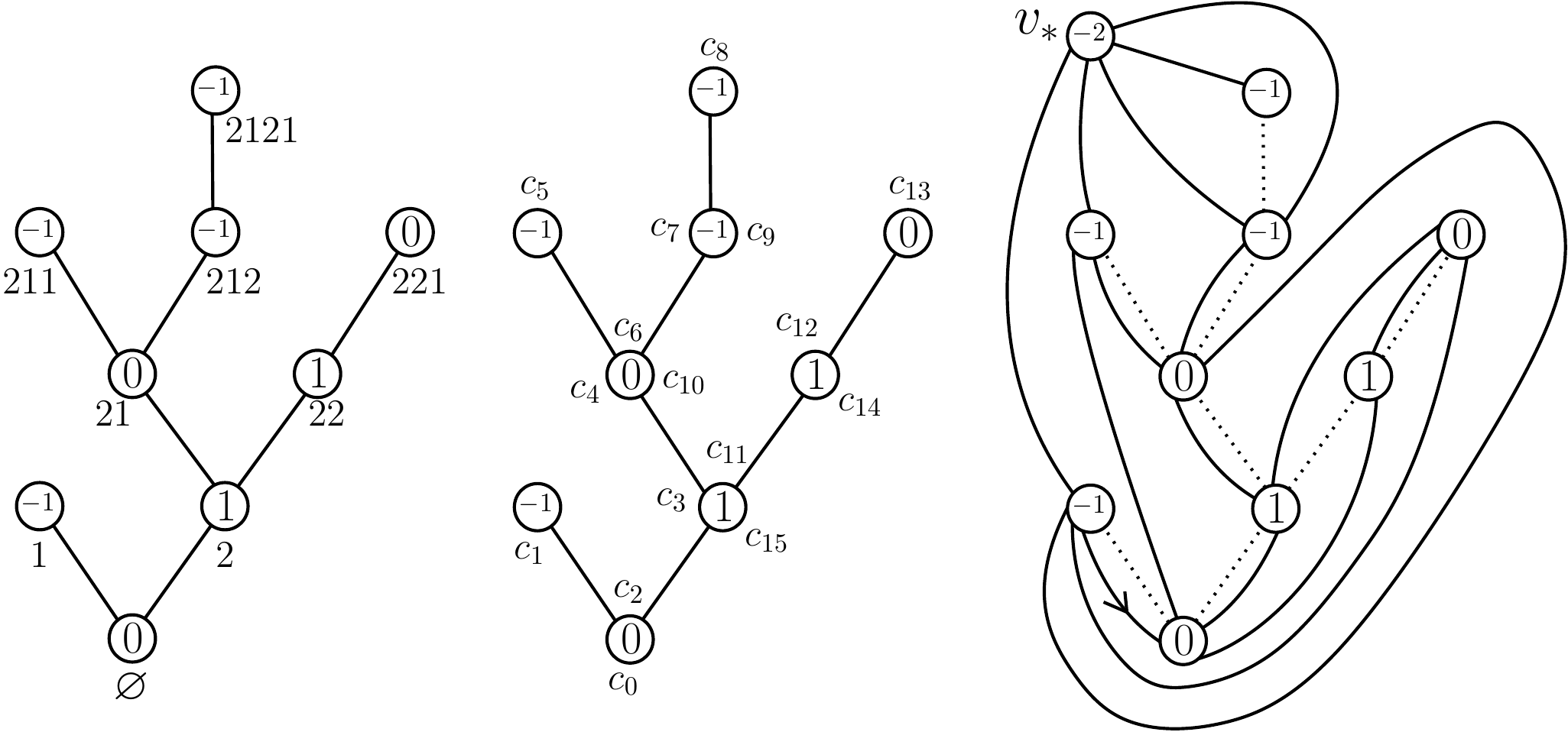}
 \caption{\label{Scha}
Schaeffer's bijection. Left: a labeled tree with $8$ edges. Middle: the sequence $c_0,c_1,\ldots,c_{15}$
of corners enumerated in cyclic order. Right: the edges of the associated quadrangulations with $8$ faces (case $\ve=-1$).}
 \end{center}
 \vspace{-5mm}
 \end{figure}

We first need to introduce the class of discrete trees that will be relevant. First recall that
a {\it plane tree} $\tau$ is a (finite) rooted ordered tree. A plane tree 
can be specified by representing each vertex as a finite word made of
positive integers, in such a way that the empty word $\varnothing$
corresponds to the root, and for instance the word $21$ corresponds to the first child
of the second child of the root. This should be clear from the left side of 
Fig.~\ref{Scha} (ignore for the moment the circled figures). To make the connection 
with planar maps, we will assume that plane trees are drawn in the plane (or rather on the sphere)
in the way illustrated in the left side of Fig.~\ref{Scha}, so that in particular the edges connecting a vertex to its parent, its first child,
its second child, etc., appear in clockwise order around that vertex. 

A {\it labeled tree} is a plane tree $\tau$, with vertex set $V(\tau)$, whose vertices are
assigned integer labels $(\ell_v)_{v\in V(\tau)}$ in such a way that the 
following two properties hold:
\begin{itemize}
\item[(i)] $\ell_\varnothing=0$;
\item[(ii)] $|\ell_v-\ell_{v'}|\leq 1$ whenever $v,v'\in V(\tau)$ are adjacent.
\end{itemize} 
The circled figures in the left side of Fig.~\ref{Scha} show a possible assignment
of labels. For every $n\geq 2$, let 
$\T_{n}$ stand for the set of all labeled trees with $n$ edges.

A rooted and pointed quadrangulation is a rooted quadrangulation given 
with a distinguished vertex (which can be any vertex, including the root vertex).
For every $n\geq 2$, let $\M^{4,\bullet}_n$ stand for the set of
all rooted and pointed quadrangulations. 

We then claim that there is a one-to-one correspondence between
the sets $\M^{4,\bullet}_n$ and $\T_{n}\times\{-1,1\}$ (this correspondence is
called Schaeffer's bijection). To explain this correspondence, let us start 
from a labeled tree $(\tau,(\ell_v)_{v\in V(\tau)})$ in $\T_{n}$
and a sign $\ve\in\{-1,+1\}$. We need to consider
corners of the tree $\tau$: A corner incident to a vertex $v$ of $\tau$
is an angular sector between two successive edges incident to $v$ (for instance,
in the tree of the left side of Fig.~\ref{Scha}, the root $\varnothing$ has $2$ corners, the vertex
$21$ has $3$ corners, and the vertex $221$ has only one corner). 
By convention, the root corner $c_0$ is the corner 
``below'' the root vertex. The set of all corners is given a cyclic ordering by
moving clockwise around the tree: starting from the root corner $c_0$, the $2n$ corners
can be listed as $c_0,c_1,\ldots,c_{2n-1}$ in cyclic ordering (see the middle part of 
Fig.~\ref{Scha}). 
We agree that every corner inherits the label of the vertex to which it
is incident.

With the labeled tree $(\tau,(\ell_v)_{v\in V(\tau)})$, we associate 
 a quadrangulation $M$ by the following device. First, the vertex set 
 of $M$ is the union of the vertex set of $\tau$ and an extra vertex $v_*$, which by
 convention is assigned the label 
 $$\ell_{v_*}=\min_{v\in V(\tau)} \ell_v - 1.$$
 Then, in order to obtain the edges of the quadrangulation $M$, we proceed in
 the following way. For every corner $c$ of $\tau$, with label $\ell_c$, we draw an edge starting from
 this corner and ending at the next corner of $\tau$ (in the cyclic ordering)
 with label $\ell_c-1$ --- this corner will be called the successor of $c$. This makes sense unless $\ell_c$ is equal to the minimal
 label on the tree $\tau$, in which case we draw an edge starting from $c$
 and ending at $v_*$. All these edges can be drawn, in a unique manner (up to homeomorphisms), in such a 
 way that they do not cross and do not cross the edges 
 of $\tau$, and the resulting planar map is a quadrangulation (see Fig.~\ref{Scha} for an example where,
 for instance, there are edges of $M$ connecting $c_0$ to $c_1$, $c_1$ to $v_*$, $c_2$ to $c_5$, $c_3$ to $c_4$, etc.).
 
 We still have to define the root of the quadrangulation and its distinguished vertex. The root edge 
 is the edge starting from $c_0$ and ending at the successor of $c_0$,
 and its orientation is determined by the sign $\ve$: The root vertex 
 is $\varnothing$ if and only if $\ve = +1$. Finally the distinguished 
 vertex of $M$ is $v_*$, and we have indeed obtained 
 a rooted and pointed quadrangulation.

 \begin{proposition}
 \label{Scha-bij}
 The preceding construction yields a bijection from $\T_{n}\times\{-1,1\}$ onto $\M^{n,\bullet}_4$.
 Moreover, if the roooted and pointed quadrangulation $M$
 is the image of the pair $((\tau,(\ell_v)_{v\in V(\tau)}),\ve)$ under this bijection,
 the vertex set $V(M)$ is canonically identified with $V(\tau)\cup \{v_*\}$ where
 $v_*$ is the distinguished vertex of $M$, and with this identification we have, for every $v\in V(\tau)$,
\begin{equation}
\label{dist-root}
\dg^M(v_*,v)= \ell_v - \min_{u\in V(\tau)} \ell_{u} +1.
\end{equation} \end{proposition}

Let us explain why property \eqref{dist-root} holds. Let $v$ be a vertex of $M$ distinct from $v_*$,
so that $v$ is identified to a vertex of $\tau$. Choose any corner $c$ incident to $v$ in the tree $\tau$. 
The construction of edges in Schaeffer's bijection shows that there is an edge connecting $c$
to a corner $c'$ of a vertex $v'$ with label $\ell_v-1$. But similarly, there is an edge of $M$ connecting the corner $c'$ of $v'$ to a 
corner of a vertex with label $\ell_v-2$.
We can continue inductively, and we get a path in $M$ of length $ \ell_v - \min_{u\in V(\tau)} \ell_{u}$  
connecting $v$ to a vertex with minimal label, which itself (by the rules of  Schaeffer's bijection)
is adjacent to $v_*$ in $M$. In this way we get the upper bound
$$\dg^M(v_*,v)\leq  \ell_v - \min_{u\in V(\tau)} \ell_{u} +1.$$
The corresponding lower bound is also very easy, using the fact that $|\ell_v-\ell_{v'}|=1$
whenever $v$ and $v'$ are adjacent in $M$, again by the construction of Schaeffer's bijection.

Property \eqref{dist-root} is useful when studying the metric properties of $M$ (in view of 
proving the case $p=4$ of Theorem \ref{main}). However, \eqref{dist-root} only gives information
about distances from the distinguished vertex $v_*$, which is far from sufficient if 
one is interested in the Gromov-Hausdorff convergence. If $v$ and $v'$ are two arbitrary vertices
of $M$, there is however a very useful upper bound for the graph distance $\dg^M(v,v')$. To
state this bound, recall that $c_0,c_1,\ldots,c_{2n-1}$ is the sequence of corners of the tree $\tau$
associated with $M$, listed in the cyclic ordering that was already used in Schaeffer's bijection.
For every $i\in\{0,1,\ldots,2n-1\}$, let $v_i$ be the vertex corresponding to the corner $v_i$.
Then, if $0\leq i<j\leq 2n-1$, we
have
\begin{equation}
\label{bound-dist}
\dg^M(v_i,v_j)\leq \ell_{v_i} + \ell_{v_j} -2 \max\Big(\min_{k\in[i,j]} \ell_{v_k}, \min_{k\in[j,2n-1]\cup[0,i]} \ell_{v_k}\Big)+2.
\end{equation}
The proof of this bound is easy. Consider the geodesic path $\gamma$ from the corner $c_i$ to $v_*$
constructed  as in the proof of \eqref{dist-root}, and the similar geodesic path from the corner $c_j$.
A simple argument shows that these two geodesic paths coalesce at a vertex whose label 
is the maximum appearing in \eqref{bound-dist} minus $1$. The concatenation of these two geodesic paths
up to their coalescence time thus gives a path from $v_i$ to $v_j$
whose length is the right-hand side of \eqref{bound-dist}.

\subsection{Ideas of the proof of Theorem \ref{main}}
\label{proof-main}

Schaeffer's bijection allows us to sketch the main ideas of the proof of Theorem \ref{main} in the case of quadrangulations. We start from
a uniformly distributed rooted and pointed quadrangulation $M_n$ with $n$ faces (the fact
that we consider a rooted and pointed quadrangulation rather than a rooted quadrangulation as
in Theorem \ref{main} is unimportant since by ``forgetting'' the distinguished vertex of $M_n$ we get a
uniformly distributed rooted quadrangulation), and we let $(\tau_n,(\ell^n_v)_{v\in V(\tau_n)})$
be the associated labeled tree. We note that $\tau_n$ is uniformly distributed over the set
of all plane trees with $n$ edges, because for every such tree there is the same number $3^n$
of possible assignments of labels. It is well known that the height of
the tree $\tau_n$ is of order $\sqrt{n}$ when $n$ is large, and, from the central limit theorem,
one may guess that the maximal and the minimal label in $\tau_n$ are of order $\sqrt{\sqrt{n}}=n^{1/4}$ (just note that 
conditionally given $\tau_n$, the increments of labels along the different edges of $\tau_n$ are independent and uniformly 
distributed over $\{-1,0,1\}$). Recalling \eqref{dist-root}, we see that the diameter of $M_n$
must be of order $n^{1/4}$, which explains the rescaling in Theorem \ref{main}. 

Then, a well-known result of Aldous shows that the tree $\tau_n$
viewed as a metric space for the graph distance rescaled by the factor $1/\sqrt{2n}$
converges in distribution to the CRT --- with our particular normalization of the CRT.
This convergence can be stated in a more precise form using the so-called
``contour functions'' which keep track of the lexicographical order on the trees. 
Furthermore, using the fact that the variance of the uniform distribution on $\{-1,0,1\}$ is $2/3$, one gets that the labels rescaled by $(2/3)^{-1/2}(2n)^{-1/4}$ converge to Brownian 
motion indexed by the CRT (we do not make the meaning of this convergence precise here).
This suggests that the scaling limit of $M_n$ can be described in terms of the
CRT equipped with Brownian labels. However, in contrast with the discrete picture, we need to
perform  some identification of vertices of the CRT. Let us explain this. Writing again $c_0,\ldots,c_{2n-1}$
for the sequence of corners of the tree $\tau_n$, we note that for $i<j$, the corner $c_i$
is connected to the corner $c_j$ by an edge of $M_n$ as soon as
$$\ell_{c_j}=\ell_{c_i}-1\hbox{ and } \ell_{c_k}\geq \ell_{c_i}\hbox{ for every }k\in\{i,i+1,\ldots,j-1\}.$$
The point is now that, even for large values of $n$, this property will hold for certain pairs $(i,j)$ such that
$j-i$ is of order $n$. Because of the rescaling of the graph distance by $n^{-1/4}$, which informally implies that
two adjacent vertices are identified in the scaling limit, this means that
certain pairs of distinct points of the CRT must be glued together.

Finally, a tightness argument relying on the bound \eqref{bound-dist} can be used to verify that sequential limits of $(V(M_n),n^{-1/4}\dg^{M_n})$
exist in the Gromov-Hausdorff sense, and are represented as quotient spaces of the CRT (equipped
with Brownian labels)
for a certain pseudo-metric $D$. The discrete bound \eqref{bound-dist} implies that
the pseudo-metric $D$ satisfies $D\leq D^\circ$, where $D^\circ$ is defined in \eqref{Dzero}.
It immediately follows that $D$ must be bounded by the right-hand side of \eqref{formulaD}.
The remaining part of the argument, which unfortunately is much harder, is to verify
that \eqref{formulaD} indeed holds. 

\section{Infinite-volume models and the Brownian plane}
\label{Bplane}

The random planar maps discussed in the preceding sections are finite (random) graphs embedded in the sphere. 
It turns out that one can also define infinite random lattices that are limits in a certain sense of 
uniformly distributed triangulations or quadrangulations with a fixed number of faces (one could consider
more general planar maps, see in particular \cite{Ste}). A pioneering work of Angel and Schramm \cite{AS},
which (together with the companion paper \cite{Ang} and the Chassaing-Schaeffer paper \cite{CS}) motivated much of the subsequent research about random planar maps,
introduced the so-called uniform infinite planar triangulation or UIPT as the local limit 
of uniformly distributed triangulations with a fixed number of faces --- in fact, Angel and Schramm
considered ``type II triangulations'' where self-loops are not allowed, but the analogous construction for general
triangulations can be found in \cite{Ste}. Let us present the analog of the Angel-Schramm construction 
for quadrangulations, which is due to Krikun \cite{Kri}.

\begin{figure}[!h]
 \begin{center}
 \includegraphics[width=10cm]{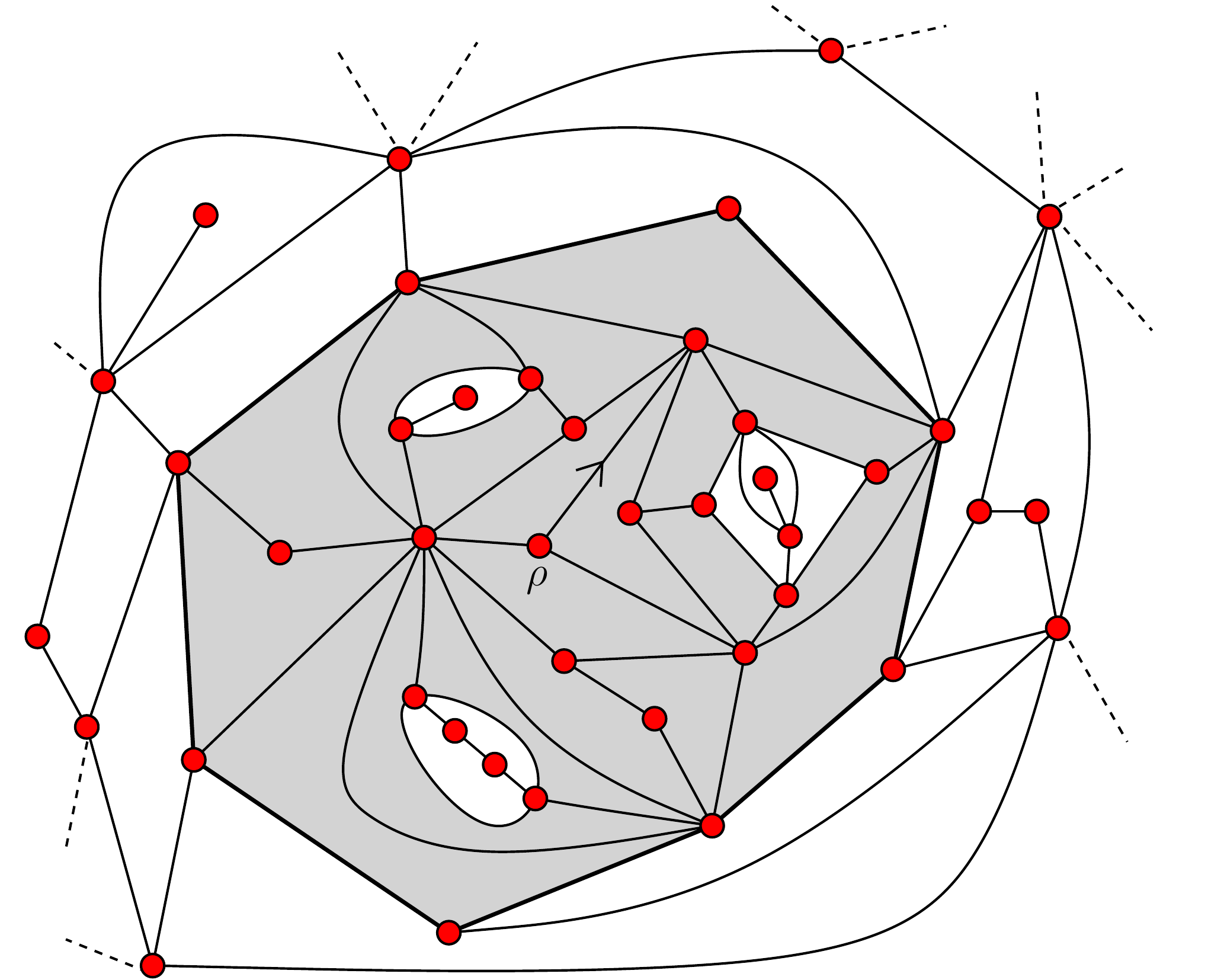}
 \caption{\label{Hull1}
A large quadrangulation $Q$ near its root vertex $\rho$ and in grey the ball $B_2(Q)$.}
\end{center}
 \vspace{-5mm}
 \end{figure}

If $M$ is a rooted planar map with root vertex $\rho$, and $r\geq 1$ is an integer, the ball 
of radius $r$ in $M$, which is denoted by $B_r(M)$, is the rooted planar map obtained by
keeping only those faces of $M$ that are incident to a vertex whose graph distance from $\rho$
is at most $r-1$. See Fig.~\ref{Hull1} for an illustration in the case of a quadrangulation. 
This definition of balls can be extended to infinite (rooted) planar lattices, meaning infinite (rooted) connected graphs properly embedded in the plane.

For every $n\geq 1$, let $Q_n$ be uniformly distributed over the set $\M^4_n$ of all rooted quadrangulations
with $n$ faces. Then one proves \cite{Kri} that there exists an infinite random rooted planar lattice $Q_\infty$
such that, for every integer $r\geq 1$ and for every rooted planar map $M$, we have
$$\P(B_r(Q_n)=M)\build{\la}_{n\to\infty}^{} \P(B_r(Q_\infty)=M).$$
The infinite random lattice $Q_\infty$ is called the 
uniform infinite planar quadrangulation or UIPQ. It is the {\it local limit} of $Q_n$ as $n\to\infty$, meaning that the distribution of what one
sees in $Q_n$ in a fixed neighborhood of the root vertex ``stabilizes'' when $n\to\infty$ to the distribution
of the corresponding neighborhood of the root vertex in $Q_\infty$. We emphasize that this convergence is very
different from the Gromov-Hausdorff convergence in Theorem \ref{main} (which also dealt with 
uniformly distributed quadrangulations): here there is no rescaling of the graph distance, and the limit
is an infinite random lattice instead of a random compact metric space. 
Both the Krikun paper \cite{Kri} and the Angel-Schramm work \cite{AS} for triangulations relied on 
enumeration techniques, but a different approach to the UIPQ based on bijections with labeled trees was proposed 
by Chassaing and Durhuus \cite{CD} (the equivalence between this approach and Krikun's one was later established
by M\'enard \cite{Men}). A simple construction of the UIPQ, relying on the version of Schaeffer's bijection 
presented in Section \ref{Scha-bi}, can be found in \cite{CMM}.

The UIPQ is an infinite-volume limit of finite quadrangulations. In the same way, one may ask about
the existence of an infinite-volume version of the Brownian map. This is the Brownian plane, which 
appears in the following theorem as a scaling limit for the UIPQ. Before stating this theorem, recall that
a metric space is called boundedly compact if every closed bounded set is compact. Write $\dg^{Q_\infty}$
for the graph distance on the vertex set $V(Q_\infty)$, and view $(V(Q_\infty),\dg^{Q_\infty})$
as a pointed metric space, where the distinguished point is the root vertex.

\begin{theorem}[\cite{Plane}]
\label{convBplane} 
There exists a random boundedly compact pointed metric space $(\mathcal{P},D_\infty)$ such that
$$(V(Q_\infty),\lambda\,\dg^{Q_\infty})\build{\la}_{\lambda\to 0}^{\mathrm{(d)}} (\mathcal{P},D_\infty),$$
where the convergence holds in distribution in the local Gromov-Hausdorff sense.
\end{theorem}

The local Gromov-Hausdorff convergence (in distribution) means that, for every real $r>0$, the closed ball of radius $r$ centered at
the distinguished point of $V(Q_\infty)$ converges (in distribution) to the same ball centered at the distinguished point in the limiting space $\mathcal{P}$, 
in the sense of the Gromov-Hausdorff distance for compact spaces. Just like the Brownian map, the Brownian plane is believed to be
a universal object, and in fact a version of the preceding theorem for the UIPT has been proved by Budzinski \cite{Bud} with the same limiting space.

We refer to \cite{Plane} for the construction of the Brownian plane, which is a continuous analog 
of the construction of the UIPQ in \cite{CMM} (a slightly different approach to the Brownian plane is given in \cite{Hull}). 
The construction of \cite{Plane} is very similar to the construction of the Brownian map in Section \ref{sec:consBM}. The key
ingredient is now Brownian motion indexed by the infinite Brownian tree, which can be understood as the Brownian tree
conditioned on non-extinction.

\begin{figure}[!h]
 \begin{center}
 \includegraphics[width=13cm]{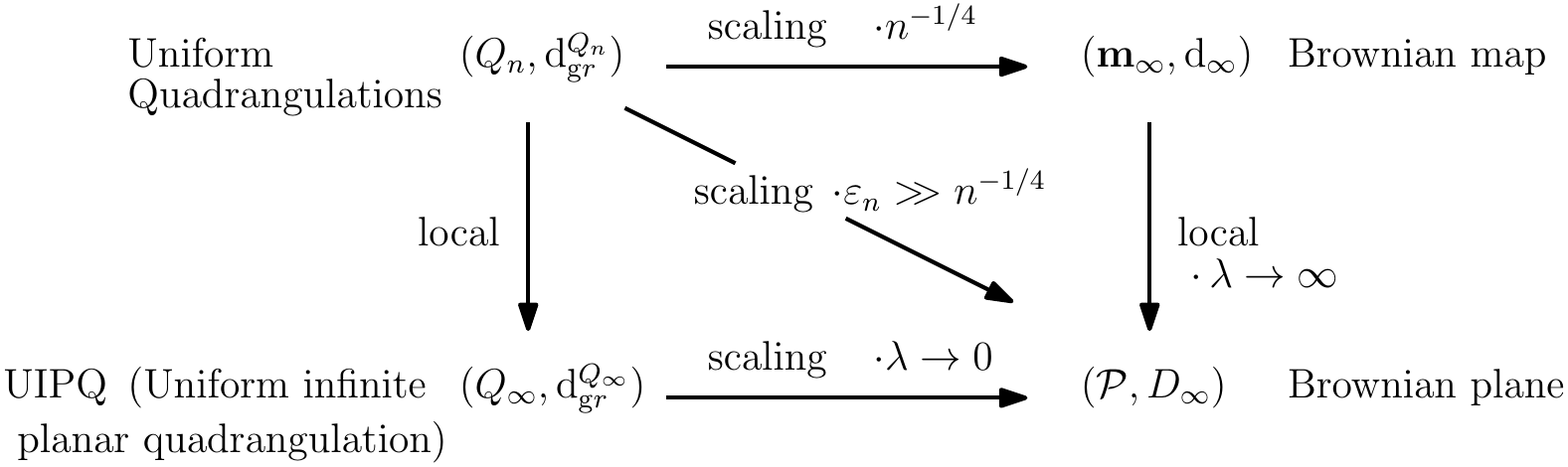}
 \caption{\label{diag}
Convergence to the Brownian plane.}
\end{center}
 \vspace{-5mm}
 \end{figure}

One may obtain the Brownian plane as a limiting object in a variety of different ways. For instance, starting from the 
Brownian map $(\bm_\infty,\dd_\infty)$ of Theorem \ref{main} and viewing $\bm_\infty$ as a pointed space with distinguished point $x_*$
(cf. the end of Section \ref{subsec:consBM}), one checks that $(\mathcal{P},D_\infty)$ is the limit
of $(\bm_\infty,\lambda\,\dd_\infty)$ when $\lambda\to \infty$, in the local Gromov-Hausdorff sense.
In the terminology of \cite{BBI}, one may say that the Brownian plane is the tangent cone (in distribution) of the Brownian map 
at $x_*$.
Alternatively one can start from the uniformly distributed quadrangulation $Q_n$ and scale the distance by
a factor $\ve_n$ tending to $0$ less fast than $n^{-1/4}$. Fig.~\ref{diag} gives a  diagram taken from
\cite{Plane} that summarizes these convergences in distribution, together with those of Theorem \ref{main} and \ref{convBplane}.

In a way similar to Theorem \ref{topo}, the Brownian plane is homeomorphic to the usual plane. On the other hand,
the Brownian plane shares the same local properties as the Brownian map (in fact in a strong sense, since one can
couple the Brownian plane and the Brownian map so that the respective balls of sufficiently small radius centered 
at the distinguished point are isometric, see \cite{Plane}). In particular, the Hausdorff dimension of the Brownian plane is also 
equal to $4$. Furthermore, the Brownian plane enjoys an additional property of scaling invariance: for every 
$\lambda>0$, the space $(\mathcal{P},\lambda\,D_\infty)$ has the same distribution as $(\mathcal{P},D_\infty)$.
This makes certain calculations more tractable in the Brownian plane than in the Brownian map: see \cite{Hull}
for several remarkable distributions related to the Brownian plane. 

\section{Planar maps with a boundary and Brownian disks}
\label{PlaBou}

In this section we introduce Brownian disks as scaling limits of
quadrangulations with a boundary. Brownian disks are
models of random geometry which unlike the Brownian map
are homeomorphic to the closed disk. Nonetheless, Brownian disks
are very closely related to the Brownian map, and, as we will discuss later, various subsets of the
Brownian map can be identified as Brownian disks.

Let us start with a basic definition. Recall that the root face of a rooted planar map
is the face lying to the left of the root edge.

\begin{definition}
\label{quad-bdry}
A quadrangulation with a (general) boundary is a rooted planar map $Q$
such that all faces but the root face have degree $4$. The root face is also called the outer face
and the other faces are called inner faces. The degree of the outer face, which is an even integer, is called the boundary size 
or the perimeter of $Q$. 
\end{definition}

\begin{figure}[!h]
 \begin{center}
 \includegraphics[width=8cm]{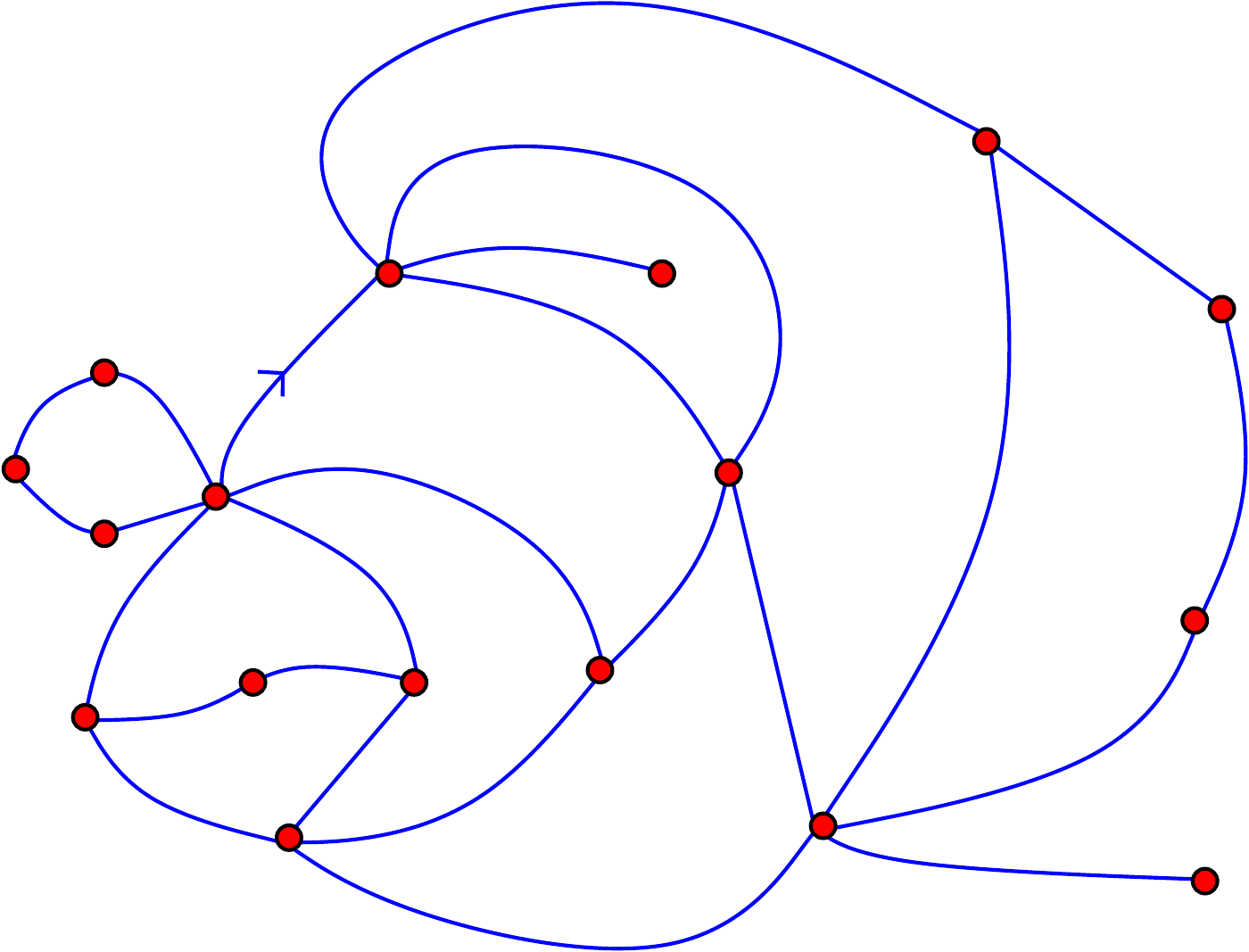}
 \caption{\label{quad-bd}
 A quadrangulation with a boundary of size 14.}
 \end{center}
 \vspace{-5mm}
 \end{figure}

See Fig.~\ref{quad-bd} for an example. One could also consider $p$-angulations
with a boundary (in particular triangulations with a boundary) but for the 
sake of simplicity we restrict our attention to quadrangulations.

For every integer $k\geq 1$, we denote the
set of all pointed quadrangulations with a boundary of size $2k$ by $\Q^{\partial,k}$.
For every integer $n\geq 0$, the subset of  $\Q^{\partial,k}$ consisting of those 
quadrangulations $Q$ that have $n$ inner faces is denoted by  $\Q^{\partial,k}_n$.
Then, for every $k\geq 1$, there is a constant $b_k>0$ such that
$$\#\Q^{\partial,k}_n \build{\sim}_{n\to\infty}^{} b_k\,12^n\, n^{-5/2}.$$
See formula (4) in \cite{CM}. 

A random variable $B_k$ with values in $\Q^{\partial,k}$
is called a Boltzmann 
quadrangulation with a boundary of size $2k$ if, for every integer $n\geq 0$
and every $Q\in \Q^{\partial,k}_n$,
$$\P(B_k=Q)= \wt b_k\,12^{-n},$$
where $\wt b_k>0$ is the appropriate normalizing constant.

The following result, which is analogous to Theorem \ref{main}, is a special case of \cite[Theorem 8]{BM}.

\begin{theorem}
\label{conv-quad-bdry}
For every integer $k\geq 1$, let $B_k$ be a Boltzmann 
quadrangulation with a boundary of size $2k$. Then,
$$\Big(V(B_k), \sqrt{3/2}\,k^{-1/2}\,\dg^{B_k}\Big) \build{\la}_{k\to\infty}^{\mathrm{(d)}} (\D, D^\partial)$$
where the convergence holds in distribution for the Gromov-Hausdorff topology. The limiting random 
compact metric space $(\D, D^\partial)$ is called the free Brownian disk with perimeter $1$.
\end{theorem}

The factor $\sqrt{3/2}$ in the convergence of the theorem is present only to allow a 
simpler description of the limit in the next section.

In contrast with Theorem \ref{main}, we notice that the number of faces of $B_k$ is not fixed, but
only its perimeter. One can verify that the number of faces of $B_k$ is typically of order $k^2$, and
so the scaling factor $k^{-1/2}$ in Theorem \ref{conv-quad-bdry} corresponds to the 
factor $n^{-1/4}$ in Theorem \ref{main}. One can prove versions of Theorem \ref{conv-quad-bdry}
for quadrangulations where both the boundary size and the volume (number of faces) are fixed
and grow to infinity simultaneously in such a way that the volume stays proportional to the square of the boundary size: This leads to the definition of Brownian disks with
given perimeter and volume. See \cite{BMR} for a discussion of all possible
scaling limits of quadrangulations with a boundary, and \cite{GM2} for an analog
of Theorem \ref{conv-quad-bdry} in the case of quadrangulations with a {\it simple} boundary.

For every $a>0$, the free Brownian disk with perimeter $a$ may be defined as
the random metric space $(\D, \sqrt{a}\,D^\partial)$.

One proves \cite{Bet} that the Brownian disk is homeomorphic to the closed unit disk, and this makes it
possible to define the boundary $\partial \D$ as the set of all points in $\D$ that have no
neighborhood homeomorphic to the open unit disk. 

\section{Excursion theory for Brownian motion indexed by the Brownian tree}
\label{excu-theory}

In this section, which is mostly taken from \cite{ALG}, we discuss an excursion theory for Brownian motion indexed by the Brownian tree.
An important motivation is to derive a construction of Brownian disks which 
is much analogous to the construction of the Brownian map explained in Section \ref{sec:consBM}.
However, we believe that this excursion theory is interesting in its own and should have
many other applications. There is of course a strong analogy with the classical It\^o theory \cite{Ito}
but also important differences due to the fact that the parameter set is a tree, and so 
connected components of the complement of the zero set of Brownian motion are $\R$-trees instead of intervals in
the classical setting. 

Recall from Definition \ref{BSexc} the $\sigma$-finite measure $\N_0(\dd \omega)$
on the space of snake trajectories with initial point $0$, and the notation 
$W_s(\omega)=\omega_s$, $\zeta_s(\omega)=\zeta_{(\omega_s)}$ for $s\geq 0$, and 
$\sigma(\omega)=\sup\{s\geq 0:\zeta_s\not =0\}$. The ``Brownian tree''
$\t_\zeta$ is the tree coded by the function $(\zeta_s)_{s\geq 0}$ as explained in Section \ref{sub:BrT},
and we use the notation $Z_a=\wh W_s$ if $a=p_\zeta(s)$, where $p_\zeta$ stands for the
canonical projection from $\R_+$ onto $\t_\zeta$. The collection $(Z_a)_{a\in \t_\zeta}$
is thus our Brownian motion indexed by the Brownian tree.

In a way very similar to classical excursion theory, our aim is
to describe the process $Z$ restricted to a connected component
of $\{b\in\t_\zeta: Z_b\not = 0\}$. To this end we first introduce the 
notion of an {\it excursion debut}. We say that $a\in \t_\zeta$ 
is an excursion debut if
\begin{itemize}
\item[(i)] $Z_a=0$;
\item[(ii)] $a$ has a strict descendant $a'$ such that $Z_b\not = 0$ for every $b\in \rrbracket a,a'\rrbracket$.
\end{itemize}
In (ii), we use the obvious notation $\rrbracket a,a'\rrbracket=\llbracket a,a'\rrbracket \backslash \{a\}$. We then observe that
connected components
of $\{b\in\t_\zeta: Z_b\not = 0\}$ are in one-to-one correspondence with excursion debuts: The connected 
component $\mathcal{C}_a$ associated with an excursion debut $a$ is just the set of all strict descendants $a'$ of $a$ such that
the property $Z_b\not = 0$ for every $b\in \rrbracket a,a'\rrbracket$ holds. 

\begin{figure}[!h]
 \begin{center}
 \includegraphics[width=11cm]{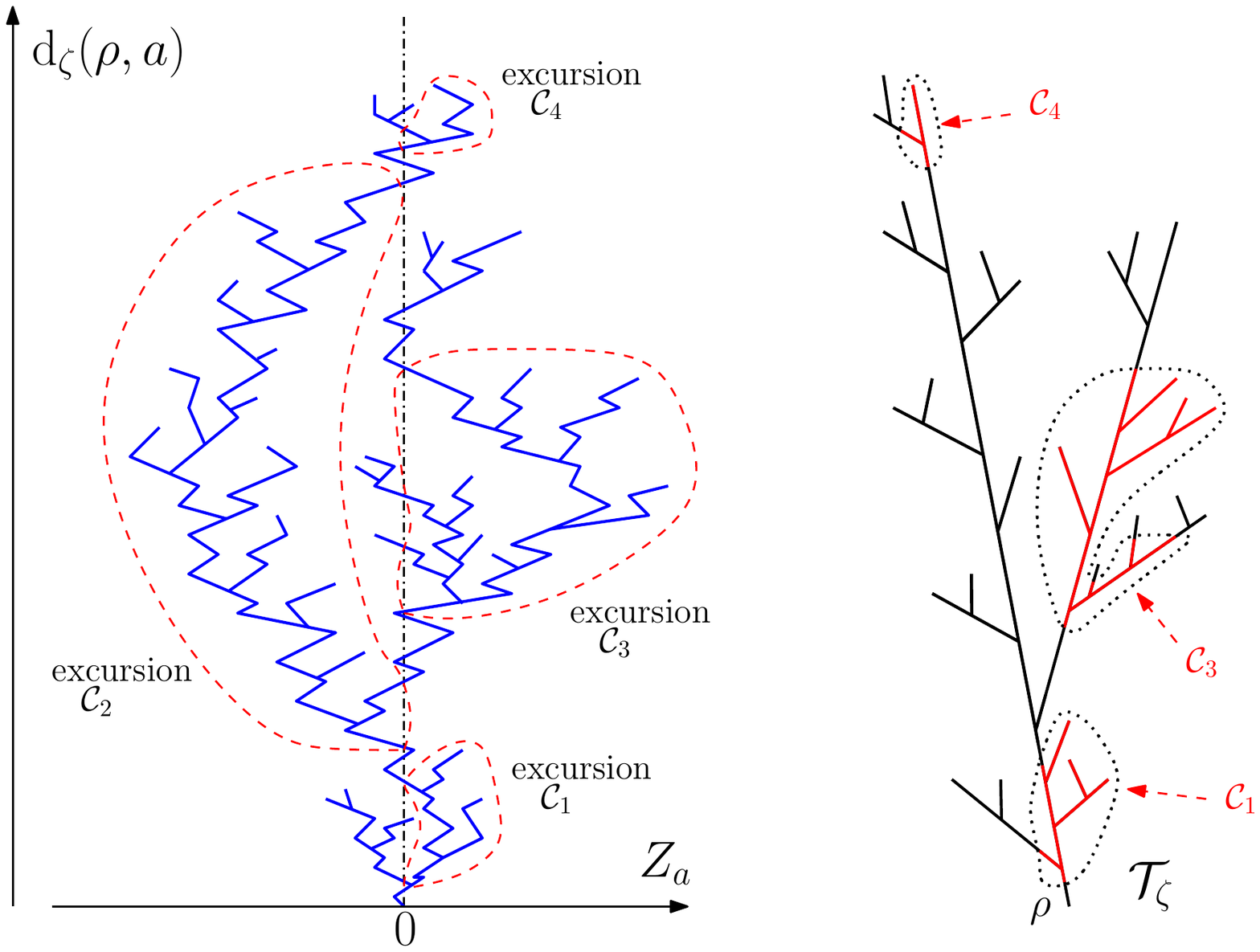}
 \caption{\label{Exc}
A schematic representation of excursions. The right side shows the tree $\t_\zeta$, and the parts of the tree
inside the dotted lines are 	a few connected components of the set $\{a\in\t_\zeta:Z_a\not =0\}$. The left side shows the
values of $Z_a$ for $a\in \t_\zeta$, or equivalently the paths $W_s$ which form a ``tree of Brownian paths'', and the parts inside the dashed lines
are a few excursions away from $0$.}
 \end{center}
 \vspace{-5mm}
 \end{figure}

We will now explain how the values of $Z$ on a given connected component can be represented by
a snake trajectory. So let us fix an excursion debut $a$. The fact that $a$ has strict descendants 
implies that there are exactly two times $u<v$ such that $p_\zeta(u)=p_\zeta(v)=a$
(there could be three such times if $a$ were a branching point of $\t_\zeta$,  but this case is excluded because branching points have nonzero labels).
Recall that $W_u=W_v$ is called the historical path of $a$. 
We note that the descendants of $a$ are exactly the points $p_\zeta(s)$ for $s\in [u,v]$.
We can then define a snake trajectory $\tilde W^{(a)}=(\tilde W^{(a)}_s)_{s\geq 0}$ in $\mathcal{S}_0$, which describes the
labels of descendants of $a$, by setting for every $s\geq 0$,
$$\tilde W^{(a)}_{s}(t):= W_{(u+s)\wedge v}(\zeta_u+t)\;,\quad \hbox{for } 0\leq t\leq \tilde\zeta^{(a)}_s:= \zeta_{(u+s)\wedge v} -\zeta_u.$$
In fact we are not interested in all descendants of $a$, but only in those that lie in
the associated connected component $\mathcal{C}_a$. For this reason, we introduce the time change
$$\eta^{(a)}_s:=\inf\Big\{r\geq 0: \int_0^r \dd t\,\mathbf{1}_{\{\tau_0^*(\tilde W^{(a)}_t)\geq \tilde\zeta^{(a)}_t\}}>s\Big\},$$
where we use the notation $\tau_0^*(\w)=\inf\{t\in(0,\zeta_{(\w)}]:\w(t)=0\}$ for $\w\in\mathcal{W}$, with the usual convention
$\inf\varnothing=+\infty$. The effect of this time change will be
to disregard the paths $\tilde W^{(a)}_s$ that return to $0$ and then survive for a positive amount of time. Setting 
for every $s\geq 0$,
$$W^{(a)}_s:=\tilde W^{(a)}_{\eta^{(a)}_s}$$
defines another snake trajectory in $\mathcal{S}_0$, which accounts for the labels 
on the connected component $\mathcal{C}_a$. We sometimes call $W^{(a)}$ the excursion
associated with the excursion debut $a$.

Let $(a_i)_{i\in I}$ be the (countable) collection of all excursion debuts. For every $i\in I$, we write 
$l_i$ for the total local time at $0$ accumulated by the historical path of $a_i$ (this makes sense
because historical paths behave like one-dimensional Brownian paths), and we note that $l_i$
is also the total local time at $0$ for the historical path of any point in the component $\mathcal{C}_{a_i}$.

\begin{theorem}[\cite{ALG}]
\label{cons-exc}
There exists a $\sigma$-finite measure $\M_0$ on $\mathcal{S}_0$ such that, for any nonnegative 
measurable function $\Phi$ on $\R_+\times \mathcal{S}_0$, we have
$$\N_0\Big(\sum_{i\in I} \Phi(l_i,W^{(a_i)})\Big) = \int _0^\infty \dd \ell\,\M_0\Big(\Phi(\ell,\cdot)\Big).$$
\end{theorem}

\begin{figure}[!h]
 \begin{center}
 \includegraphics[width=9cm]{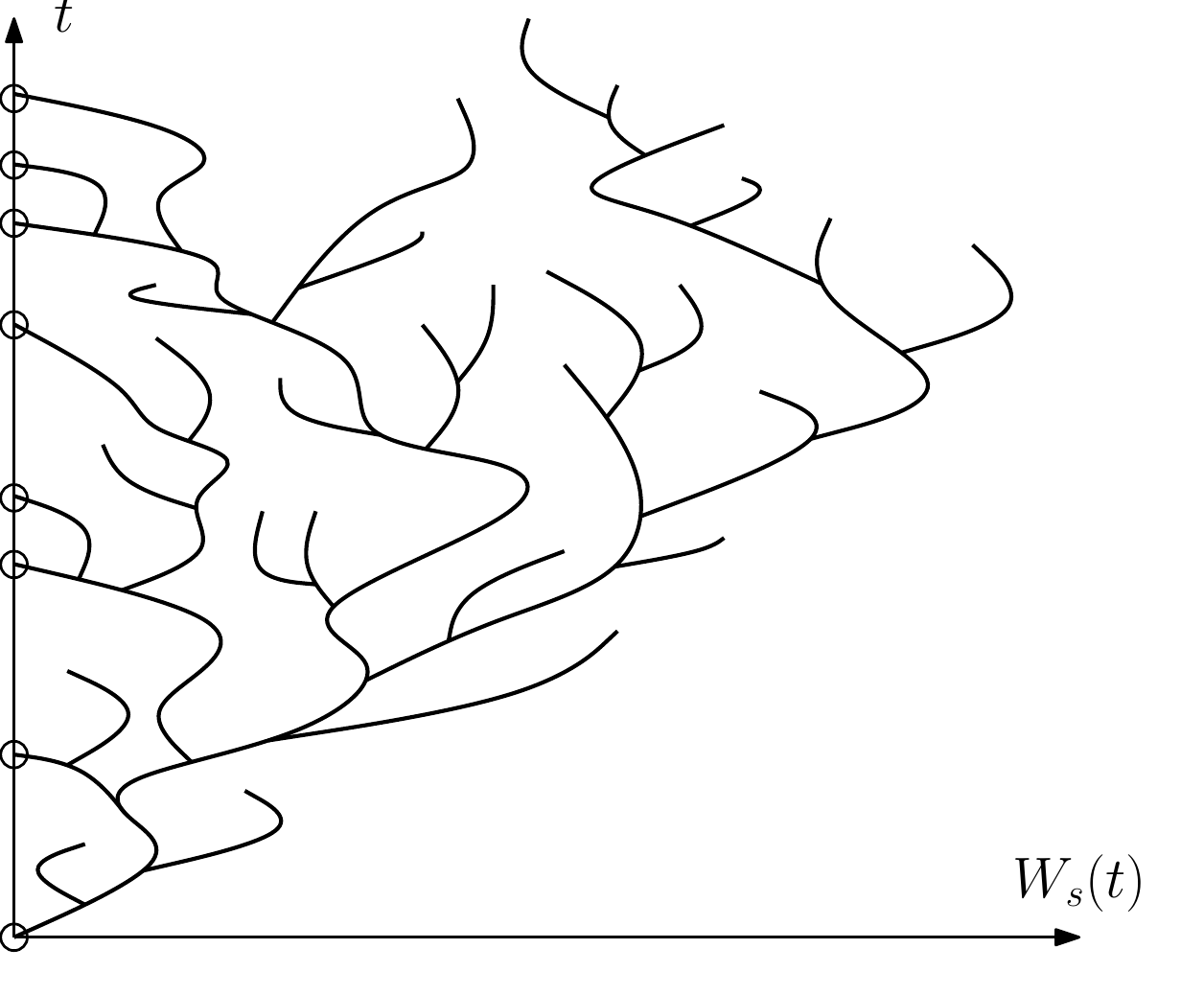}
 \caption{\label{treepaths}
A schematic representation of the paths $W_s$ under $\N^*_0$. The quantity $\mathcal{Z}^*_0$ measures the ``number''
of circled points corresponding to returns of certain paths $W_s$ to $0$.}
 \end{center}
 \vspace{-5mm}
 \end{figure}

For symmetry reasons, we may write
$$\M_0= \frac{1}{2}\Big( \N_0^* + \check \N^*_0\Big)$$
where $\N^*_0$ is supported on nonnegative snake trajectories, and $\check \N^*_0$ is the push forward of 
$\N^*_0$ under the mapping $\omega\mapsto -\omega$. Under $\N^*_0$,
the paths $W_s$ form a ``tree of Brownian paths'' starting from $0$, which take positive values until the first time
when they return to $0$ (if they do return to $0$) and are stopped at that time if not earlier. See Fig.~\ref{treepaths}
for a schematic illustration. 

Theorem \ref{cons-exc} provides a first-moment formula for the collection of excursions $(W^{(a_i)})_{i\in I}$, but,
in contrast with the classical excursion theory, this result does not say anything about the independence of these excursions.
To discuss independence properties, we first need to introduce the ``boundary size'' of an excursion, which roughly speaking 
measures the quantity of paths $W_s$ that return to $0$.

\begin{proposition}
\label{bdry-size}
The limit 
$$\mathcal{Z}^*_0:=\lim_{\ve\to 0} \frac{1}{\ve^{2}}\int_0^\sigma \mathbf{1}_{\{0<|\wh W_s|<\ve\}}\,\dd s$$
exists $\M_0$ a.e.
\end{proposition}

Using scaling arguments, it is not hard to define the conditional probability measures 
$\M_0(\cdot \mid \mathcal{Z}^*_0=z)$ for every $z>0$. 

In order to state the main result of this section, we still need to introduce a process $(\Lambda_r)_{r> 0}$ 
defined under the excursion measure $\N_0$, such that, for every $r>0$, $\Lambda_r$ ``counts the number'' of
paths $W_s$ that accumulate a total local time $r$ at $0$ and are stopped when they have accumulated that amount of
local time. The precise definition of $\Lambda_r$ fits in the general framework of exit measures as presented in
\cite[Chapter V]{Zurich}, but can also be given via the following approximation:
$$\Lambda_r:=\lim_{\ve\to 0} \frac{1}{\ve} \int_0^\sigma \mathbf{1}_{\{\chi_r(W_s)<\zeta_s<\chi_r(W_s)+\ve\}} \dd s\,,\qquad \N_0\ \hbox{a.e.}$$
where $\chi_r(W_s)=\inf\{t\geq 0: L^0_t(W_s)>r\}$, if $(L^0_t(W_s))_{0\leq t\leq \zeta_s}$ denotes the local time at $0$
of the path $W_s$. 

Thanks to the special Markov property of the Brownian snake (see the appendix of \cite{subor}), one can prove that the process $(\Lambda_r)_{r>0}$
is Markovian under $\N_0$ (this makes sense even though $\N_0$
is an infinite measure because $\N_0(\Lambda_r\not =0)<\infty$ for every $r>0$) with the transition kernels
of the continuous-state branching process with branching mechanism $\psi(u)=\sqrt{8/3}\,u^{3/2}$. In
particular, $(\Lambda_r)_{r>0}$ has a c\`adl\`ag modification with only positive jumps, which we
consider in the next statement.

Recall that $l_i$ denotes the total local time at $0$ accumulated by the historical path of $a_i$.

\begin{theorem}[\cite{ALG}]
\label{excu-theo}
The numbers $l_i$, $i\in I$, are exactly the jump times of the process $(\Lambda_r)_{r>0}$. Furthermore,
conditionally on the process $(\Lambda_r)_{r>0}$, the excursions $W^{(a_i)}$, $i\in I$, are independent
and, for every $j\in I$, the conditional distribution of $W^{(a_j)}$ is $\M_0(\cdot\mid \mathcal{Z}^*_0 =\Delta \Lambda_{l_j})$.
\end{theorem}

In particular, the boundary size of the excursion $W^{(a_j)}$ is $\Delta \Lambda_{l_j}$.

In the applications developed below, we will be interested mainly in positive excursions and in the measure $\N^*_0$,
which we call the positive Brownian snake excursion measure. As in the case of $\M_0$
we can define the conditional probability measures
$$\N^{*,z}_0:= \N^*_0(\cdot\mid \mathcal{Z}^*_0=z)$$
for every $z>0$. 

Interestingly, a number of explicit distributions can be computed explicitly under $\N^*_0$.
In particular the joint distribution of the pair $(\mathcal{Z}^*_0,\sigma)$ (boundary size and volume) under $\N^*_0$ 
has a density
given by
$$f(z,s)=\frac{\sqrt{3}}{2\pi} \, \sqrt{z}\,s^{-5/2}\,\exp\Big(-\frac{z^2}{2s}\Big).$$
Consequently, for every fixed $z>0$, the density of $\sigma$ under $\N^{*,z}_0$ is 
$$g_z(s)= \frac{1}{\sqrt{2\pi}}\,z^3\,s^{-5/2}\,\exp\Big(-\frac{z^2}{2s}\Big).$$
The latter density also appears as the density of the asymptotic distribution of the
rescaled volume (number of faces) of a Boltzmann quadrangulation with 
perimeter $2k$, when $k\to\infty$. This will be explained by the results of the
next section. 

In the classical setting of excursions away from $0$ for a standard linear Brownian motion
starting from $0$, it
is well known that the process can be reconstructed by concatenating the 
different excursions (some care is required since there are infinitely many
excursions on any interval $[0,t]$, $t>0$). In our setting of a tree-indexed process, things are more complicated 
since excursions 
are no longer ordered linearly, but have a certain genealogical structure induced by the genealogy of their debuts: 
In the example of Fig.~\ref{Exc}, the excursion $\cc_1$ is an ancestor of both $\cc_3$ and $\cc_4$, but $\cc_3$ is not ancestor of $\cc_4$. 
Still this genealogical structure can be described in the following way. 

For every $a,a'\in\t_\zeta$, we let $\delta(a,a')$ be the total local time at $0$ accumulated by
the process $Z$ along the line segment $\llbracket a,a'\rrbracket$ of the tree $\t_\zeta$. This makes sense since 
we know that $Z$ evolves like a linear Brownian motion along any segment of the tree. Then $\delta(\cdot,\cdot)$ is
a pseudo-metric on $\t_\zeta$, and we can define the associated equivalence 
relation by setting $a\simeq a'$
if and only if $\delta(a,a')=0$. Obviously $a\simeq a'$ holds if $a$ and $a'$ belong to the
same connected component of $\{b\in\t_\zeta: Z_b\not = 0\}$ (because then $Z$
does not vanish on $\llbracket a,a'\rrbracket$). 
The quotient space $\t_\zeta/\!\simeq$ can thus be seen as obtained from $\t_\zeta$
by gluing each excursion component into a single point. It turns out \cite{subor} that this quotient space
has a remarkable probabilistic structure.

\begin{theorem}[\cite{subor}]
\label{subo-tree}
Under $\N_0$, the quotient space $\t_\zeta/\!\simeq$ equipped with the distance induced by $\delta$ is a stable L\'evy tree with index $3/2$.
\end{theorem}

\begin{figure}[!h]
 \begin{center}
 \includegraphics[width=11cm]{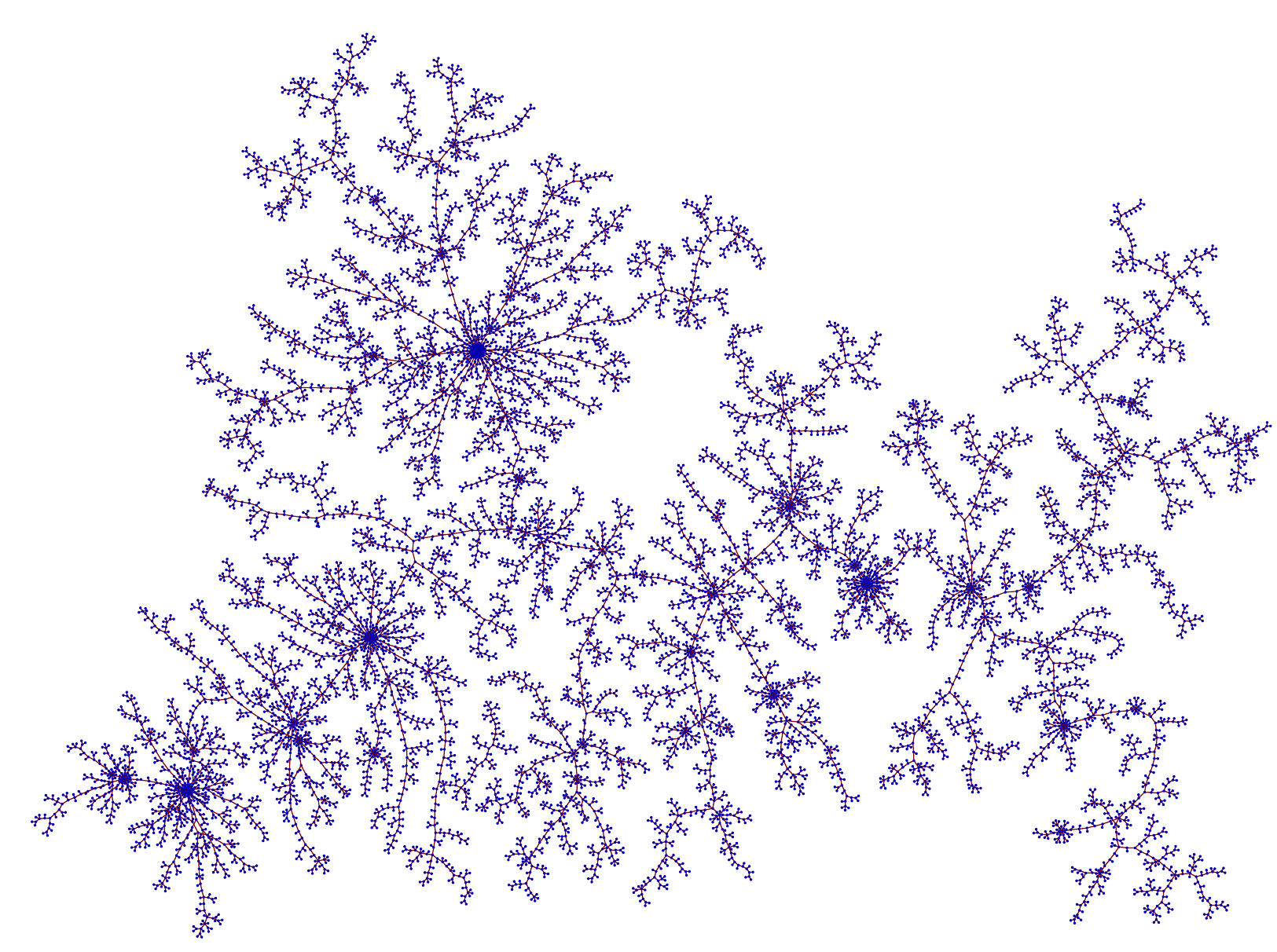}
 \caption{\label{stable-tree}
A simulation of the stable tree with index $3/2$ (simulation: I. Kortchemski).}
 \end{center}
 \vspace{-5mm}
 \end{figure}

We refer to \cite{DLG} for the definition and main properties of L\'evy trees (note that
these trees are typically defined under an infinite measure). See 
Fig.~\ref{stable-tree} for a simulation. The stable L\'evy tree
appearing in the theorem has infinitely many points of infinite multiplicity, and to each such
point one can assign a ``mass'' corresponding informally to the degree of the
point in the tree. Then one can check that 
points of infinite multiplicity of the tree $\t_\zeta/\!\simeq$ are in one-to-one correspondence with excursions
of $(Z_a)_{a\in\t_\zeta}$, and that the mass of every point of infinite multiplicity coincides with the
boundary size of the associated excursion.

The previous lines then suggest the following possible ``reconstruction'' method
(which we will not attempt to make rigorous here). Starting from 
a stable L\'evy tree with index $3/2$, associate independently with each point $a$ of infinite 
multiplicity a snake trajectory $W^a$ distributed according to $\M_0(\cdot\mid \mathcal{Z}^*_0=m_a)$, where 
$m_a$ is the mass of $a$, then ``glue'' at the location of  $a$ the genealogical tree of the
snake trajectory $W^a$, with the corresponding labels inherited from $W^a$. The resulting random $\R$-tree
equipped with labels should be the Brownian tree equipped with Brownian labels.

\section{Constructing Brownian disks from the positive Brownian snake excursion measure}
\label{Br-Disk}

The results of this section are taken from \cite{Disk}. 
The first naive idea to construct a free Brownian disk is to imitate the construction
of Section \ref{subsec:consBM}, replacing the measure $\N_0$
by $\N^*_0$. This does not give the desired result, but yields another
object of interest, namely the (free) Brownian disk with glued boundary.

Recall our notation $(\D,D^\partial)$ for the free Brownian disk with perimeter $1$,
and $\partial \D$ for the boundary of $\D$. We define a pseudo-metric on $\D$
by setting, for every $x,y\in\D$,
$$D^\dagger(x,y)=\min\{D^\partial(x,y), D^\partial(x,\partial \D)+D^\partial(y,\partial \D)\}.$$
Clearly $D^\dagger(x,y)=0$ if and only if $x=y$, or both $x$ and $y$ belong to
$\partial \D$. Write $\D^\dagger$ for the set obtained from $\D$ by identifying 
all points of the boundary $\partial \D$ to a single point. Then $D^\dagger$
induces a metric on $\D^\dagger$, which we still denote by $D^\dagger$. The
compact metric space $(\D^\dagger,D^\dagger)$ is called the free Brownian disk
with perimeter $1$ and glued boundary. The case of a perimeter equal to $z$
is treated analogously.

\begin{proposition}
\label{cons-glued}
The random metric space $(\mm,D)$ defined
via the construction of Section \ref{subsec:consBM} from a 
snake trajectory $\omega$ distributed according
to $\N^{*,z}_0$ is a free Brownian disk
with perimeter $z$ and glued boundary.
\end{proposition}

The problem is then to recover the free Brownian disk from the
same object with glued boundary. This can indeed be achieved by
a slight modification of the construction of Section \ref{subsec:consBM}. 

From now on, we argue under the measure $\N^{*,z}_0(\dd\omega)$
for some fixed $z>0$. Recalling that $\t_\zeta$
denotes the genealogical tree of the snake trajectory $\omega$, we use the same
notation $Z_a$ for the ``labels'' on $\t_\zeta$ ($Z_a=\wh W_s$ if $a=p_\zeta(s)$). In contrast 
with the case of $\N_0$, labels are now nonnegative reals, and we define
the ``boundary'' $\partial \t_\zeta$ by
$$\partial \t_\zeta:=\{a\in\t_\zeta: Z_a=0\}.$$

Recalling the definition of $D^\circ$ in \eqref{Dzero}, we set, for every $a,b\in \t_\zeta\backslash \partial \t_\zeta$,
$$\Delta^\circ(a,b)=\left\{
\begin{array}{ll}
D^\circ(a,b)\qquad&\hbox{if } \displaystyle{\max\Big(\min_{c\in[a,b]} Z_c,\min_{c\in[b,a]}Z_c\Big)>0,}\\
\infty&\hbox{otherwise}.
\end{array}\right.$$
Roughly speaking, the condition in the first line of the last display means that we can go
from $a$ to $b$ ``around'' the tree $\t_\zeta$ without visiting a vertex of $\partial\t_\zeta$. 
We then define $\Delta(a,b)$ for every $a,b\in \t_\zeta\backslash \partial \t_\zeta$ by the exact analog
of formula \eqref{formulaD}:
$$
\Delta(a,b) = \inf\Big\{ \sum_{i=1}^k \Delta^\circ(a_{i-1},a_i)\Big\},
$$
where the infimum is over all choices of the integer $k\geq 1$ and of the
elements $a_0,a_1,\ldots,a_k$ of $\t_\zeta\backslash \partial \t_\zeta$ such that $a_0=a$
and $a_k=b$. 
One easily verifies that
the mapping $(a,b)\mapsto \Delta(a,b)$ takes finite values and is continuous on 
$(\t_\zeta\backslash \partial \t_\zeta)\times (\t_\zeta\backslash \partial \t_\zeta)$.

\begin{theorem}
\label{BrDisk}
With probability one under $\N^{*,z}_0$, the function $(a,b)\mapsto \Delta(a,b)$
has a continuous extension to $\t_\zeta\times \t_\zeta$, which is a pseudo-metric on $\t_\zeta$.
We let $\Theta$ stand for the associated quotient space, and we equip $\Theta$ with the induced
metric, which is still denoted by $\Delta(a,b)$. Then, the random metric space 
$(\Theta,\Delta)$ is a free Brownian disk with perimeter $z$ under $\N^{*,z}_0$,
and its boundary $\partial\Theta$ is the image of $\partial\t_\zeta$
under the canonical projection from $\t_\zeta$ onto $\Theta$. Furthermore, 
if $x\in\Theta$ is the image of $a\in\t_\zeta$ under the canonical projection, we have
$$\Delta(x,\partial \Theta)=Z_a.$$
\end{theorem}

We note that we can define a 
volume measure $\mathbf{V}(\dd x)$ on $\Theta$ as the image of the volume measure on $\t_\zeta$ under
the canonical projection. In particular the total mass of $\mathbf{V}$ is $\mathbf{V}(\Theta)=\sigma$ (recall our notation
$\sigma$ for the duration of the snake trajectory $\omega$, which is also the total
mass of the volume measure on $\t_\zeta$). Hence we may define the 
Brownian disk with perimeter $z$ and volume $v$ as the random metric space 
$(\Theta,\Delta)$ under the conditional probability measure $\N^{*,z}_0(\cdot \mid \sigma =v)$.
This is consistent with the construction of \cite{Bet,BM} using scaling limits of
quadrangulations with a boundary with fixed perimeter and volume.

A nice feature of the construction of Theorem \ref{BrDisk} (in contrast with the previous constructions in \cite{Bet,BM})
is the fact that labels $Z_a$ now correspond to distances from the boundary. This also makes it
possible to construct a natural ``length measure'' on the boundary. 
The following proposition is closely related to the approximation of $\mathcal{Z}^*_0$ in Proposition \ref{bdry-size}.

\begin{proposition}
\label{measure-boundary}
Almost surely under $\N^{*,z}_0$, there exists a finite measure $\nu$
on $\partial\Theta$ with total mass $z$, such that, for every bounded
continuous function $\varphi$ on $\Theta$,
$$\langle \nu, \varphi\rangle = \lim_{\ve\to 0} \frac{1}{\ve^2}\int_\Theta \bV(\mathrm{d}x)\,\varphi(x)\,\mathbf{1}_{\{\Delta(x,\partial\Theta)<\ve\}}.$$
\end{proposition}

We will now exhibit certain particular subsets of the Brownian map that are Brownian disks. 
So we now argue under the measure $\N^{(1)}_0(\dd \omega)$ and consider the
metric space $(\mm,D)$ constructed in Section \ref{subsec:consBM}. Recall 
from the end of this section that $\mm$ has a distinguished point $x_*$ such that
distances from $x_*$ exactly correspond to the labels $Z_x$ up to a shift
(see \eqref{dist-min} above). We will discuss properties of the connected components 
of the complement of balls centered at $x_*$. At this point, we should mention that
the point $x_*$ does not play a special role, and that the re-rooting invariance properties
of the Brownian map \cite[Section 8]{Geo} show that the same properties
hold if $x_*$ is replaced by a point chosen according to the volume measure on
the Brownian map. We recall that this volume measure, which is denoted by $\mathbf{v}(\dd x)$,
is the push forward of the volume measure on $\t_\zeta$, and that
$\mathbf{v}$ is a probability measure under $\N^{(1)}_0(\dd \omega)$.

We note that the Brownian map is a length space (as a Gromov-Hausdorff limit of length spaces)
and that, if $O$ is an open subset of $\mm$,
we can define an {\it intrinsic} metric $D^{O}_{\mathrm{intr}}$ on $O$ by declaring that
$D^{O}_{\mathrm{intr}}(x,y)$ is the minimal length of a continuous path connecting $x$ to $y$ in $O$
(see \cite[Chapter 2]{BBI}). 

For every $z>0$ and $v>0$, we let $\F_{z,v}$ be the distribution of the Brownian disk with perimeter
$z$ and volume $v$. The following statement can be found in \cite[Theorem 3]{Disk} (see also \cite{JM}
for a related work).

\begin{theorem}
\label{ccBm}
Let 
$r>0$ and let $B(x_*,r)$ stand for the closed ball of radius $r$ centered at $x_*$
in $(\mm,D)$. Then, $\N^{(1)}_0$ a.s. for every connected component $\mathbf{C}$ of $\mm\backslash B(x_*,r)$, the 
limit
\begin{equation}
\label{approx-bd}
|\partial\mathbf{C}|:=\lim_{\ve\to 0} \frac{1}{\ve^2} \int_\mathbf{C} \bv(\mathrm{d}x)\,\mathbf{1}_{\{ D(x,\partial\mathbf{C})<\ve\}}
\end{equation}
exists and is called the boundary size of $\mathbf{C}$. On the event $\{\mm\backslash B(x_*,r)\not =\varnothing\}$, write $\mathbf{C}^{r,1},\mathbf{C}^{r,2},\ldots$
for the connected components of $\mm\backslash B(x_*,r)$ ranked in decreasing order of their boundary sizes.
Let $D^{r,j}_{\mathrm{intr}}$ be the intrinsic distance on $\mathbf{C}^{r,j}$.
Then, $\N^{(1)}_0$ a.s. on the event $\{\mm\backslash B(x_*,r)\not =\varnothing\}$, for every $j=1,2,\ldots$,
the metric $D^{r,j}_{\mathrm{intr}}$ has a continuous extension to the closure
$\ov{\mathbf{C}}^{r,j}$ of $\mathbf{C}^{r,j}$, and this extension is a metric on $\ov{\mathbf{C}}^{r,j}$. Furthermore, 
under $\N^{(1)}_0(\cdot\mid \mm\backslash B(x_*,r)\not =\varnothing)$ and conditionally on  the sequence 
$$(|\partial \mathbf{C}^{r,1}|,\bv(\mathbf{C}^{r,1})), (|\partial \mathbf{C}^{r,2}|,\bv(\mathbf{C}^{r,2})),\ldots$$
the metric spaces $(\ov{\mathbf{C}}^{r,j},D^{r,j}_{\mathrm{intr}})$, $j=1,2,\ldots$, are independent
Brownian disks with respective distributions $\F_{|\partial \mathbf{C}^{r,j}|,\bv(\mathbf{C}^{r,j})}$, $j=1,2,\ldots$.
\end{theorem}

Let us briefly explain why Theorem \ref{ccBm} is related to the excursion theory developed in
Section \ref{excu-theory}. The key point is the fact that distances from $x_*$ are given 
(up to the  shift by $-Z_*$) by the labels $Z_x$. Assuming that $r>-Z_*$ for simplicity, it is then not too hard to verify that connected components
of the complement of $B(x_*,r)$ correspond --- via the construction presented in Section \ref{sec:consBM} --- to excursions of Brownian motion indexed by
the Brownian tree above the (random) level $r+Z_*$. The distribution of these excursions
can be analysed thanks to Theorem \ref{excu-theo} and we also use Theorem \ref{BrDisk}
to relate the positive Brownian snake excursion measure to the law of Brownian disks. There
are however two significant technical difficulties, because on one hand we have to deal with
excursions above a {\it random} level, instead of level $0$ in Theorem \ref{excu-theo},
and on the other hand, we argue under $\N^{(1)}_0$ instead of $\N_0$ in Section \ref{excu-theory}.

Informally, Theorem \ref{ccBm} says that connected components of the complement of 
a ball centered at a ``typical point'' in the Brownian map are independent Brownian disks
conditionally on their boundary sizes and volumes. A similar result \cite[Theorem 18]{Disk} holds 
for the connected components of the complement of the Brownian net, which is a particular subset 
of the free Brownian map playing an important role in the axiomatic characterization 
of Miller and Sheffield \cite{MS0}. At this point, we mention that we could have stated a version 
of Theorem \ref{ccBm} for the free Brownian map, which is nicer in the sense that we
do not need to condition on the volumes: We get that the connected components of the
complement of a ball centered at $x_*$ are independent free Brownian disks conditionally
on their perimeters. In the next section, we discuss a similar statement for 
the free Brownian disk, where distances from $x_*$ are replaced by distances
from the boundary.

\section{Slicing Brownian disks at heights}
\label{CutBr}

In this section, which is based on  \cite{LGR}, we consider the random metric space $(\Theta,\Delta)$
defined in Theorem \ref{BrDisk}, which is a free Brownian disk
with perimeter $z$ under the probability measure $\N^{*,z}_0$. 
For every $x\in \Theta$, define the {\it height} of $x$ by
$$H(x)=\Delta(x,\partial \D).$$
We also consider the maximal height
$$H^*= \max_{x\in \Theta} H(x).$$
Recall the notation $\bV(\dd x)$ for the volume measure on $\Theta$. 

\begin{theorem}
\label{excu}
Let $r>0$. Then, $\N^{*,z}_0$ a.s., for every connected component 
$\cc$ of $\{x\in\Theta:H(x)>r\}$, the limit
$$|\partial \cc|=\lim_{\ve\to 0} \frac{1}{\ve^2} \int_\mathcal{C} \bV(\mathrm{d}x)\,\mathbf{1}_{\{ H(x)<r+\ve\}}$$
exists and is called the perimeter of $\cc$.
On the event $\{H^*>r\}$, let $\cc^{r,1},\cc^{r,2},\ldots$ be the connected components of 
$\{x\in\Theta:H(x)>r\}$ ranked in decreasing order of their perimeters. Then, a.s.
on the event $\{H^*>r\}$, for every $j=1,2,\ldots$, the intrinsic metric on $\cc^{r,j}$ has
a continuous extension to the closure
$\ov{\mathcal{C}}^{r,j}$ of $\mathcal{C}^{r,j}$, which is a metric on $\ov{\mathcal{C}}^{r,j}$, and 
conditionally on the perimeters $|\partial \cc^{1,r}|,|\partial\cc^{2,r}|,\ldots$,
the resulting metric spaces $\ov{\cc}^{r,1},\ov{\cc}^{r,2},\ldots$ are independent free Brownian disks.
\end{theorem}

As explained for Theorem \ref{ccBm} at the end of the previous section, Theorem \ref{excu}
can be derived from the excursion theory developed in Section \ref{excu-theory}. The difficulty 
now comes from the fact that we must argue under $\N^{*,z}_0$ instead of $\N_0$
in Section \ref{excu-theory}.

With the notation of Theorem \ref{excu},  an obvious question is to describe the distribution of the process 
$$\mathbf{X}(r)= (|\partial \cc^{1,r}|,|\partial\cc^{2,r}|,\ldots)$$
giving for every $r>0$ the perimeters of all connected components of $\{x\in\Theta:H(x)>r\}$ (by
convention $\mathbf{X}(r)=(0,0,\ldots)$ if $H^*\leq r$). 
We also take $\mathbf{X}(0)=(z,0,0,\dots)$ and view $(\mathbf{X}(r))_{r\geq 0}$ as
a random process taking values in the space of nonincreasing sequences of 
nonnegative real numbers. 
Theorem \ref{excu} then suggests that this process enjoys properties
similar to those of the growth-fragmentation processes that have been studied
recently by several authors. In fact, Bertoin, Curien and Kortchemski \cite{BCK} (see
also \cite{BBCK} for extensions) have considered a process analogous to $\mathbf{X}$
for triangulations with a boundary and showed that the scaling limit
of this process (when the boundary size tends to infinity) is a
well-identified growth-fragmentation process. Still it does not seem easy
to apply the results of \cite{BCK} in order to identify the distribution of the process $\mathbf{X}$,
but the excursion theory of Section \ref{excu-theory} can be used instead to compute this
distribution.

Before stating our last result, we need to recall a few basic facts about growth-fragmentation
processes (see \cite{Ber} for more details). The starting ingredient is a positive self-similar 
Markov process $(Y_t)_{t\geq 0}$ with only negative jumps. Suppose that $Y_0=z$, and view $(Y_t)_{t\geq 0}$
as the evolution in time of the mass of an initial particle also called the Eve particle. 
At each time $t$ where the process $Y$ has a jump, we consider that a new particle with mass $-\Delta Y_t$ (a
child of the Eve particle)
is born, and the mass of this new particle evolves (from time $t$) again according to the
law of the process $Y$, but independently of the evolution of the mass of the Eve particle. Then 
each child of the Eve particle has children at discontinuity times of its mass process, and so on.
Under suitable assumptions (see \cite{Ber}), we can make sense of the process $(\mathbf{Y}(t))_{t\geq 0}$ 
giving for every time $t$ the sequence  (in
decreasing order) of masses of all particles alive at that time. The process  $\mathbf{Y}$ 
is Markovian and is called the growth-fragmentation
process with Eve particle process $Y$.

\begin{theorem}
\label{GF}
Under $\N^{*,z}_0$, the process $(\mathbf{X}(r))_{r\geq 0}$ is a growth-fragmentation process, which is constructed 
from an Eve particle process $X$
whose distribution starting from $1$ is specified as follows:
$$X_t=\exp(\xi_{\tau(t)}),$$
where 
$$\tau(t)=\inf\Big\{u\geq 0: \int_0^u e^{\xi_s/2}\,\dd s >t\Big\}$$
and $\xi$ is the spectrally negative L\'evy process
such that, for every $q>0$, $\E[\exp(q\,\xi_t)]=\exp(t \psi(q))$, with
\begin{equation}
\label{formula-psi}
\psi(q)= \sqrt{\frac{3}{2\pi}}\,\Bigg(
- \frac{8}{3}\,q + \int_{1/2}^1 (x^q-1+q(1-x))\,(x(1-x))^{-5/2}\,\dd x\Bigg).
\end{equation}
\end{theorem}

\noindent{\bf Remark.} The process $\xi$ drifts to $-\infty$ and the event $\{\tau(t)=\infty\}$
occurs with positive probability if $t>0$: on this event, we of course make the convention that $\exp(\xi_{\infty})=0$. 

\medskip

The expression of the process $X$ in terms of the L\'evy process $\xi$ is a special case of the classical 
Lamperti representation of a positive self-similar 
Markov process (here with index $1/2$) in terms of a
spectrally negative L\'evy
process. 
The formula for $\psi$ is the same as formula (1) in \cite{BCK}, except for the (unimportant)
multiplicative constant $\sqrt{\frac{3}{2\pi}}$. This should not come as a surprise in
view of preceding comments.

Although we have chosen to state them as properties of the free Brownian disk, Theorems
\ref{excu} and \ref{GF} are really results about the
tree-indexed Brownian motion $(Z_a)_{a\in\t_\zeta}$ under $\N^{*,z}_0$. In particular,
Theorem \ref{GF} relies on the identification of the distribution of the process giving, for each $r\geq 0$, the 
sequence of boundary sizes of all excursions above level $r$ of  the process 
$(Z_a)_{a\in\t_\zeta}$ under $\N^{*,z}_0$. There is a striking analogy with the fragmentation process
occuring when cutting the CRT at a fixed height: Precisely, it is shown in \cite{Ber0} that the sequence
of volumes 
of the connected components of the complement of the ball of radius $r$ centered at the root in the CRT
is a self-similar fragmentation process whose dislocation measure has
the form $(2\pi)^{-1/2} (x(1-x))^{-3/2}\,\dd x$, to be compared with the 
measure $(x(1-x))^{-5/2}\,\dd x$ appearing in formula \eqref{formula-psi}. 

As a consequence of Theorem \ref{GF} and known asymptotics \cite[Corollary 4.5]{BBCK} for the distribution of the
extinction time of a growth-fragmentation process, we derive the following 
corollary about the tail of the distribution of the maximal height in a Brownian disk.

\begin{corollary}
\label{max-height}
There exist positive constants $c_1$ and $c_2$ such that, for every $r\geq 1$,
$$c_1\,r^{-6}\leq \N^{*,z}_0(H^*>r)\leq c_2\,r^{-6}.$$
\end{corollary}

\end{document}